\documentclass[10pt,a4paper]{amsart}
\usepackage{amssymb,amsmath}
\usepackage[american,french]{babel}
\usepackage{amsopn}
\usepackage{graphicx}
\usepackage{fancyhdr}
\usepackage[T1]{fontenc}
\title{About quadratic differential operators}
\newcommand{\rr}{\mathbb{R}}
\newcommand{\eps}{\varepsilon}
\newcommand{\nn}{\mathbb{N}}
\newcommand{\cc}{\mathbb{C}}

\newcommand{\lde}{L^2(\rr^n)}

\def\init{\setcounter{equa}{0}}
\def\inc{\stepcounter{equa}}
\def\num{\tag{\thesection.\theequa}}

\def\wrtext#1{\relax\ifmmode{\leavevmode\hbox{#1}}\else{#1}\fi}

\def\begeq{\begin{equation}}
\def\endeq{\end{equation}}

\def\part#1{\frac{\partial}{\partial #1}}

\def\norm#1{||\,#1\,||}

\def\og{orthogonal}

\begin{document}
\newcounter{equa}
\selectlanguage{american}
\begin{center}
{\large\textbf{SUBELLIPTIC ESTIMATES FOR QUADRATIC DIFFERENTIAL OPERATORS}\\
\bigskip
\medskip
Karel \textsc{Pravda-Starov}\\
\bigskip
Imperial College, London}
\end{center}
\bigskip
\bigskip

\newtheorem{lemma}{Lemma}[subsection]
\newtheorem{definition}{Definition}[subsection]
\newtheorem{proposition}{Proposition}[subsection]
\newtheorem{theorem}{Theorem}[subsection]

\textbf{Abstract.} We prove global subelliptic estimates for quadratic differential operators. Quadratic differential operators are operators defined in the Weyl quantization by complex-valued quadratic symbols. In a previous joint work with M.~Hitrik,  
we pointed out the existence of a particular linear subvector space in the phase space intrinsically associated to their Weyl symbols, 
called singular space, which rules spectral properties of non-elliptic quadratic operators. The purpose of the present paper is to prove that quadratic operators whose singular spaces are 
reduced to zero, are subelliptic with a loss of \og derivatives \fg \  depending directly on particular algebraic properties of the Hamilton maps of their Weyl symbols.
More generally, when singular spaces are symplectic spaces, we prove that quadratic operators are subelliptic in any direction of the symplectic orthogonal complements of their singular spaces.

\medskip

\noindent
\textbf{Key words.} Quadratic differential operators, subelliptic estimates, singular space, Wick quantization. 

\medskip

\noindent
\textbf{2000 AMS Subject Classification.} 35B65, 35S05.

\section{Introduction}
\init

\subsection{Miscellaneous facts about quadratic differential operators}

Since the classical work by J. Sj\"ostrand~\cite{sjostrand}, the study of spectral properties of quadratic diffe\-rential
operators has played a basic r\^ole in the analysis of partial differential operators with double characteristics.
Roughly speaking, if we have, say, a classical pseudodifferential operator $p(x,\xi)^w$ on $\rr^n$ with the Weyl
symbol 
$p(x,\xi)=p_m(x,\xi)+p_{m-1}(x,\xi)+\ldots$ of order $m$, and if $X_0=(x_0,\xi_0)\in \rr^{2n}$ is a point
where 
$$p_m(X_0)=dp_m(X_0)=0,$$ 
then it is natural to consider the quadratic form $q$ which begins the Taylor expansion of
$p_m$ at $X_0$ in order to investigate the properties of the pseudodifferential operator $p(x,\xi)^w$. For example, the study of a priori estimates such as hypoelliptic estimates of the form
$$
\norm{u}_{m-1}\leq C_K \left(\norm{p(x,\xi)^w u}_0+\norm{u}_{m-2}\right),\quad u\in C^{\infty}_0(K), \quad
K\subset \subset \rr^n,
$$
then often depends on the spectral analysis of the quadratic operator $q(x,\xi)^w$.
See also~\cite{hypoelliptic}, as well as Chapter 22 of~\cite{hormander} together with further references given there.
In~\cite{sjostrand}, the spectrum of a general quadratic differential operator has been determined, under the
basic assumption of global ellipticity for the associated quadratic form.

\medskip
In a recent joint work with M.~Hitrik, we investigated spectral properties of non-elliptic quadratic operators. 
Quadratic operators are pseudodifferential operators
defined in the Weyl quantization,
\begin{equation}\label{3}\inc
q(x,\xi)^w u(x) =\frac{1}{(2\pi)^n}\int_{\rr^{2n}}{e^{i(x-y).\xi}q\Big(\frac{x+y}{2},\xi\Big)u(y)dyd\xi}, \num
\end{equation}
by some symbols $q(x,\xi)$, where $(x,\xi) \in \rr^{n} \times \rr^n$ and $n \in \nn^*$, which
are complex-valued quadratic forms. Since these symbols are quadratic forms, the corresponding operators in
(\ref{3}) are in fact differential operators. Indeed, the Weyl quantization of the quadratic symbol
$x^{\alpha} \xi^{\beta}$, with $(\alpha,\beta) \in \nn^{2n}$ and $|\alpha+\beta| \leq 2$, is the differential operator
$$
\frac{x^{\alpha}D_x^{\beta}+D_x^{\beta} x^{\alpha}}{2}, \ D_x=i^{-1}\partial_x.
$$
One can also notice that quadratic differential operators
are a priori formally non-selfadjoint since their Weyl symbols in (\ref{3}) are complex-valued.

\medskip

Considering quadratic operators whose Weyl symbols have real parts with a sign, say here, Weyl symbols with non-negative real parts 
\begin{equation}\label{smm1}\inc
\textrm{Re }q \geq 0, \num
\end{equation}
we pointed out in \cite{hps} the existence of a particular linear subvector space $S$ in the phase space $\rr_x^n \times \rr_{\xi}^n$ intrinsically associated to their Weyl symbols $q(x,\xi)$, 
called singular space, which seems to play a basic r\^ole in the understanding of the properties of these non-elliptic quadratic operators.
We first proved in \cite{hps} (Theorem~1.2.1) that when the singular space $S$
has a symplectic structure then the associated heat equation
\begin{equation}\label{sm1}\inc
\left\lbrace
\begin{array}{c}
\displaystyle \frac{\partial u}{\partial t}(t,x)+q(x,\xi)^w u(t,x)=0  \\
u(t,\textrm{\textperiodcentered})|_{t=0}=u_0 \in L^2(\rr^n),
\end{array} \right.
 \num
\end{equation}
is smoothing in every direction of the orthogonal complement $S^{\sigma \perp}$ of $S$ with respect to the canonical
symplectic form $\sigma$ on $\rr^{2n}$,
\begin{equation}
\label{11}\inc
\sigma \big{(}(x,\xi),(y,\eta) \big{)}=\xi.y-x.\eta, \ (x,\xi) \in \rr^{2n},  (y,\eta) \in \rr^{2n}, \num
\end{equation}
that is, that, if $(x',\xi')$ are some linear symplectic coordinates on the symplectic space $S^{\sigma \perp}$
then we have for all $t>0$, $N \in \nn$ and $u \in L^2(\rr^n)$,
\begin{equation}\label{sm4bis}\inc
\big((1+|x'|^2+|\xi'|^2)^N\big)^we^{-t q(x,\xi)^w}u \in L^2(\rr^n). \num
\end{equation}
We also proved in \cite{hps} (Theorem~1.2.2) that when the Weyl symbol $q$ of a quadratic operator fulfills (\ref{smm1}) and an assumption of
partial ellipticity on its singular space $S$ in the sense that
\begin{equation}\label{sm2}\inc
(x,\xi) \in S, \ q(x,\xi)=0 \Rightarrow (x,\xi)=0, \num
\end{equation}
then this singular space always has a symplectic structure and the spectrum of the operator
$q(x,\xi)^w$ is only composed of a countable number of eigenvalues of finite multiplicity, with a structure
similar to the one known in the case of global ellipticity~\cite{sjostrand}.

\medskip

In the present paper, we are interested in investigating the r\^ole played by the singular space when studying subelliptic properties of quadratic operators. We shall first prove that 
quadratic operators whose singular spaces are reduced to zero, fulfill global subelliptic estimates 
\begin{equation}\label{lol1}\inc
\big\|\big(\langle(x,\xi)\rangle^{2(1-\delta)}\big)^w u\big\|_{L^2} \lesssim \|q(x,\xi)^w u\|_{L^2}+\|u\|_{L^2},\num
\end{equation}
where $\langle(x,\xi)\rangle=(1+|x|^2+|\xi|^2)^{1/2}$, with a loss of \og derivatives\fg \  $\delta>0$ which  
can be directly characterized by algebraic conditions on the Hamilton maps of their Weyl symbols. More generally, when singular spaces $S$ have a symplectic structure, we prove that quadratic operators are subelliptic in any direction of the symplectic orthogonal complements of their singular spaces $S^{\sigma \perp}$, in sense that, if $(x',\xi')$ are some linear symplectic coordinates on $S^{\sigma \perp}$ then
\begin{equation}\label{lol1b}\inc
\big\|\big(\langle(x',\xi')\rangle^{2(1-\delta')}\big)^w u\big\|_{L^2} \lesssim \|q(x,\xi)^w u\|_{L^2}+\|u\|_{L^2},\num
\end{equation}
where again, the loss of \og derivatives\fg \  $\delta'>0$  
can be directly characterized by algebraic conditions on the Hamilton maps of their Weyl symbols.

\medskip
Before giving the precise statement of our main result, we shall recall miscellaneous facts and notations about quadratic differential operators. In all the following, we consider
\begin{eqnarray*}
q : \rr_x^n \times \rr_{\xi}^n &\rightarrow& \cc\\
 (x,\xi) & \mapsto & q(x,\xi),
\end{eqnarray*}
a complex-valued quadratic form with a non-negative real part
\begin{equation}\label{inf1}\inc
\textrm{Re }q(x,\xi) \geq 0, \ (x,\xi) \in \rr^{2n}, n \in \nn^*. \num
\end{equation}
We know from \cite{mehler} (p.425) that the
maximal closed realization of  the operator $q(x,\xi)^w$, i.e., the operator on $L^2(\rr^n)$ with the domain
$$D(q)=\big\{u \in L^2(\rr^n) :  q(x,\xi)^w u \in L^2(\rr^n)\big\},$$
coincides with the graph closure of its restriction to $\mathcal{S}(\rr^n)$,
$$q(x,\xi)^w : \mathcal{S}(\rr^n) \rightarrow \mathcal{S}(\rr^n).$$
Associated to the quadratic symbol $q$ is the numerical range $\Sigma(q)$ defined as the closure in the
complex plane of all its values,
\begin{equation}\label{9}\inc
\Sigma(q)=\overline{q(\rr_x^n \times \rr_{\xi}^n)}. \num
\end{equation}
We also recall from~\cite{hormander} that the Hamilton map $F \in M_{2n}(\cc)$ associated to the quadratic form $q$
is the map uniquely defined by the identity
\begin{equation}\label{10}\inc
q\big{(}(x,\xi);(y,\eta) \big{)}=\sigma \big{(}(x,\xi),F(y,\eta) \big{)}, \ (x,\xi) \in \rr^{2n},  (y,\eta) \in \rr^{2n}, \num
\end{equation}
where $q\big{(}\textrm{\textperiodcentered};\textrm{\textperiodcentered} \big{)}$ stands for the polarized form
associated to the quadratic form $q$. It directly follows from the definition of the Hamilton map $F$ that
its real part $\textrm{Re } F$ and its imaginary part $\textrm{Im }F$ are the Hamilton maps associated
to the quadratic forms $\textrm{Re } q$ and $\textrm{Im }q$, respectively. One can notice from (\ref{10}) that a
Hamilton map is always skew-symmetric with respect to $\sigma$. This is just a consequence of the
properties of skew-symmetry of the symplectic form and symmetry of the polarized form,
\begin{equation}\label{12}\inc
\forall X,Y \in \rr^{2n}, \ \sigma(X,FY)=q(X;Y)=q(Y;X)=\sigma(Y,FX)=-\sigma(FX,Y).\num
\end{equation}

Associated to the symbol $q$, we defined in~\cite{hps} its singular space $S$ as
the following intersection of kernels,
\begin{equation}\label{h1}\inc
S=\Big(\bigcap_{j=0}^{+\infty}\textrm{Ker}\big[\textrm{Re }F(\textrm{Im }F)^j \big]\Big) \cap \rr^{2n}, \num
\end{equation}
where the notations $\textrm{Re } F$ and $\textrm{Im }F$ stand respectively for the real part and the imaginary part of the Hamilton map associated to $q$.
Notice that the Cayley-Hamilton theorem applied to $\textrm{Im }F$ shows that
$$(\textrm{Im }F)^k X \in \textrm{Vect}\big(X,...,(\textrm{Im }F)^{2n-1}X\big), \ X \in \rr^{2n}, \ k \in \nn,$$
where $\textrm{Vect}\big(X,...,(\textrm{Im }F)^{2n-1}X\big)$ is the vector space spanned by the vectors $X$, ...,
$(\textrm{Im }F)^{2n-1}X$, and therefore
the singular space is actually equal to the following finite intersection of the kernels,
\begin{equation}\label{h1bis}\inc
S=\Big(\bigcap_{j=0}^{2n-1}\textrm{Ker}\big[\textrm{Re }F(\textrm{Im }F)^j \big]\Big) \cap \rr^{2n}. \num
\end{equation}

\subsection{Statement of the main results}

In this paper, we shall first study the specific case where the singular space $S$ is reduced to $\{0\}$. By assuming that 
\begin{equation}\label{h1biskrist1}\inc
S=\{0\}, \num
\end{equation}
we can therefore consider the smallest integer $0 \leq k_0 \leq 2n-1$ such that
\begin{equation}\label{h1biskrist2}\inc
\Big(\bigcap_{j=0}^{k_0}\textrm{Ker}\big[\textrm{Re }F(\textrm{Im }F)^j \big]\Big) \cap \rr^{2n}=\{0\}, \num
\end{equation}
and state the following result:

\bigskip

\begin{theorem}\label{theorem1}
Consider a quadratic operator $q(x,\xi)^w$ whose Weyl symbol  
\begin{eqnarray*}
q : \rr_x^n \times \rr_{\xi}^n &\rightarrow& \cc\\
 (x,\xi) & \mapsto & q(x,\xi),
\end{eqnarray*}
is a complex-valued quadratic form fulfilling \emph{(\ref{inf1})} and \emph{(\ref{h1biskrist1})} then the operator $q(x,\xi)^w$ fulfills the following global subelliptic estimate 
\begin{equation}\label{kristen1}\inc
\exists C>0, \forall u \in D(q), \  \big\|\big(\langle(x,\xi)\rangle^{2/(2k_0+1)}\big)^w u\big\|_{L^2} \leq C\big(\|q(x,\xi)^w u\|_{L^2}+\|u\|_{L^2}\big), \num
\end{equation}
where  $k_0$ stands for the smallest integer $0 \leq k_0 \leq 2n-1$ such that \emph{(\ref{h1biskrist2})} is fulfilled, and $\langle(x,\xi)\rangle=(1+|x|^2+|\xi|^2)^{1/2}$.
\end{theorem}

\bigskip

We shall begin our few comments about the result of Theorem~\ref{theorem1} by noticing that the estimate (\ref{kristen1}) is easy to obtain in the case where $k_0=0$. Indeed, we shall check in the following, that in this case the operator $q(x,\xi)^w$ is necessarily elliptic, and we recall from~\cite{sjostrand} that, when  
$q(x,\xi)^w$ is an elliptic quadratic operator whose Weyl symbol fulfill (\ref{inf1}),\footnote{One can actually only assume that $\Sigma(q) \neq \cc$, when $n=1$, see Lemma~3.1 in~\cite{sjostrand}.}
that is, an operator whose Weyl symbol $q(x,\xi)$ is globally elliptic on the phase space $\rr^{2n}$,
\begin{equation}\label{sm3}\inc
(x,\xi) \in \rr^{2n}, \ q(x,\xi)=0 \Rightarrow (x,\xi)=0, \num
\end{equation}
then one can construct a parametrixe inducing that this elliptic quadratic operator defines a Fredholm operator of index 0 with discrete spectrum (Theorem 3.5 in \cite{sjostrand}),
\begin{equation}\label{14}\inc
q(x,\xi)^w : B \rightarrow \lde, \num
\end{equation}
where $B$ is the Hilbert space 
\inc
\begin{align*}\label{14.1}
B= & \ \big\{ u \in \lde : q(x,\xi)^wu \in \lde \big\} \num \\
= & \ \big\{ u \in \lde : x^{\alpha} D_x^{\beta} u \in \lde \ \textrm{if} \ |\alpha+\beta| \leq 2\big\}, 
\end{align*}
with the norm  
$$\|u\|_B^2=\sum_{|\alpha+\beta| \leq 2}{\|x^{\alpha} D_x^{\beta} u\|_{\lde}^2}.$$
We therefore have in this case the natural a priori estimate 
\begin{equation}\label{lay30}\inc
\exists C>0, \forall u \in B, \  \big\|\big(\langle(x,\xi)\rangle^{2}\big)^w u\big\|_{L^2} \leq C\big(\|q(x,\xi)^w u\|_{L^2}+\|u\|_{L^2}\big),\num
\end{equation}
which gives the estimate (\ref{kristen1}) of Theorem~\ref{theorem1} when $k_0=0$.

A noticeable example of quadratic operator fulfilling the assumptions of Theorem~\ref{theorem1} is the  
Fokker-Planck operator
$$K=-\Delta_v+\frac{v^2}{4}-\frac{1}{2}+v.\partial_x-\big(\partial_xV(x)\big).\partial_v, \ (x,v) \in \rr^{2},$$
with a quadratic potential
$$V(x)=\frac{1}{2}ax^2, \ a \in \rr^*.$$
Here, we consider this non-elliptic operator only in the one-dimensional case, but it is of course just for convenience reasons.
Considering this example, our Theorem~\ref{theorem1} allows to recover the global subelliptic estimate proved by B.~Helffer and F.~Nier in \cite{HN} (Proposition~5.22),
\begin{equation}\label{kristen2}\inc
\exists C>0, \forall u \in D(K), \  \|\Lambda_x^{2/3} u\|_{L^2}^2+\|\Lambda_v u\|_{L^2}^2 \leq C\big(\|K u\|_{L^2}^2+\|u\|_{L^2}^2\big), \num
\end{equation}
where 
$$\Lambda_x=(-\Delta_x+x^2/4)^{1/2} \textrm{ and } \Lambda_v=(-\Delta_v+v^2/4)^{1/2}.$$
The Fokker-Planck operator with a quadratic potential can indeed be expressed as
$$K=q(x,v,\xi,\eta)^w-\frac{1}{2},$$
with a Weyl symbol
$$q(x,v,\xi,\eta)=\eta^2+\frac{1}{4}v^2+i(v \xi-a x \eta),$$
which is a non-elliptic complex-valued quadratic form whose real part is non-negative.
By checking that the associated Hamilton map 
$$q(x,v,\xi,\eta)=\sigma\big((x,v,\xi,\eta),F(x,v,\xi,\eta) \big),$$
is given by 
$$F= \left( \begin{array}{cccc}
  0 & \frac{1}{2}i & 0 & 0 \\
  -\frac{1}{2}ai& 0 & 0 & 1 \\
0 & 0 & 0 &  \frac{1}{2}ai \\
0 & -\frac{1}{4} & -\frac{1}{2}i & 0 
  \end{array}
\right),$$
and that the singular space 
$$S=\textrm{Ker}(\textrm{Re }F) \cap \textrm{Ker}(\textrm{Re }F \ \textrm{Im }F) \cap \rr^{4},$$
is equal to $\{0\}$, we therefore deduce from Theorem~\ref{theorem1} the global subelliptic estimate 
\begin{equation}\label{kristen3}\inc
\exists C>0, \forall u \in D(K), \  \big\|\big(\langle(x,v,\xi,\eta)\rangle^{2/3}\big)^w u\big\|_{L^2} \leq C(\|K u\|_{L^2}+\|u\|_{L^2}). \num
\end{equation}
Notice that the improvement in the variables $(v,\eta)$ appearing in the estimate (\ref{kristen2}) is easily obtained by using the Cauchy-Schwarz inequality in the following estimate
$$2\|\Lambda_v u\|_{L^2}^2-\|u\|_{L^2}^2=2\textrm{Re}(Ku,u) \leq 2\|Ku\|_{L^2}\|u\|_{L^2} \leq \|Ku\|_{L^2}^2 + \|u\|_{L^2}^2.$$

The work of B.~Helffer and F.~Nier in~\cite{HN} about this particular example of the Fokker-Planck operator with a quadratic potential  
has been the starting point of our investigation of subelliptic properties for quadratic differential operators. Nevertheless, the reader will notice that our proof of Theorem~\ref{theorem1} 
will not use the same approach as the one followed by B.~Helffer and F.~Nier. Indeed,  the proof of (\ref{kristen2}) in \cite{HN} really takes advantage of the very specific structure of the Fokker-Planck operator and seems difficult to adapt in a general setting. For our proof, we shall rather use a multiplier method inspired from the work of F.~H\'erau, J.~Sj\"ostrand and C.~Stolk in \cite{HeSjSt}, once we will have achieved the construction of a weight function (Proposition~\ref{prop1}).

One can explain the loss of \og derivatives \fg \ (See (\ref{lol1})), $\delta=2k_0/(2k_0+1)$ appearing in the estimate (\ref{kristen1}) by the following informal discussion. There are two different types of points $X_0=(x_0,\xi_0)$ in the phase space $\rr^{2n}$: those for which $\textrm{Re }q(X_0)>0$ and those for which $\textrm{Re }q(X_0)=0$. Difficulties will come from the presence of this second type of points, and the fact that the set
$\partial \Sigma(q) \cap \Sigma_{\infty}(q)$,
where (See Theorem~1.4 in~\cite{DeSjZw}),
$$ \Sigma_{\infty}(q)=\Big\{z \in \cc : z=\lim_{j \rightarrow +\infty} q(x_j,\xi_j), \ |(x_j,\xi_j)| \rightarrow +\infty \textrm{ when } j \rightarrow +\infty\Big\},$$
may not be empty in general.
In order to deal with that kind of points, we shall take advantage from the noticeable property that the average of the real part of $q$,
\begin{equation}\label{lol2}\inc
\langle \textrm{Re }q \rangle_T(X)=\frac{1}{2T}\int_{-T}^T{\textrm{Re }q(e^{tH_{\textrm{Im}q}}X)dt} \gg |X|^2, \num
\end{equation}
by the flow generated by the Hamilton vector field of its imaginary part
$$H_{\textrm{Im}q}=\frac{\partial \textrm{Im }q}{\partial \xi}.\frac{\partial}{\partial x}-\frac{\partial \textrm{Im }q}{\partial x}.\frac{\partial}{\partial \xi},$$
is always a positive definite quadratic form when its singular space $S=0$. This particular property (proved in \cite{hps}) ensures that the operator $q(x,\xi)^w$ is of principal-type
$$d \textrm{Im }q(X_0) \neq 0,$$
in any non-zero point $X_0 \in \rr^{2n}$ for which $\textrm{Re }q(X_0)=0$. We also noticed in~\cite{hps} (See Remark, Section~2) that the property (\ref{lol2}) induces that one can find for any 
non-zero point $X_0 \in \rr^{2n}$ such that $\textrm{Re }q(X_0)=0$, a positive integer $1 \leq k \leq 2n-1$ such that 
\begin{equation}\label{lol3}\inc
\forall \ 0 \leq j \leq 2k-1, \ H_{\textrm{Im}q}^{j}\textrm{Re }q(X_0)=0 \textrm{ and } H_{\textrm{Im}q}^{2k}\textrm{Re }q(X_0)\neq0.\num
\end{equation}
All the points are therefore of finite type. Since moreover the condition $(P)$ holds because of the sign property of $\textrm{Re }q$, one can microlocalize 
the operator $q(x,\xi)^w$ in a neighborhood of a point $X_0 \in \rr^{2n}$ such that (\ref{lol3}) holds, to the subelliptic model operator with large parameter $\Lambda \geq 1$,
$$D_t+i\Lambda^2 t^{2k},$$
where roughly speaking, $\Lambda \sim (x^2+\xi^2)^{1/2}$; for which the classical a priori estimate
\begin{equation}\label{lol4}\inc
\|D_tu+i\Lambda^2 t^{2k}u\|_{L^2} \gtrsim (\Lambda^2)^{\frac{1}{2k+1}}\|u\|_{L^2}, \num
\end{equation}
is fulfilled.  This informal discussion allows to understand from (\ref{lol4}) from where the loss of \og derivatives \fg \ appearing in (\ref{kristen1}) comes.  Indeed, the integer $k_0$ in Theorem~\ref{theorem1} that we characterize there by other algebraic properties on the Hamilton map, can also be characterized as the smallest integer $0 \leq k_0 \leq 2n-1$ such that for any $X \in \rr^{2n}$,  $X \neq 0$, 
$$\exists \ 0\leq k\leq k_0, \forall \ 0 \leq j \leq 2k-1, \ H_{\textrm{Im}q}^{j}\textrm{Re }q(X)=0 \textrm{ and } H_{\textrm{Im}q}^{2k}\textrm{Re }q(X)\neq0.$$ 

\bigskip

Let us now consider the more general case where the singular space $S$ defined in (\ref{h1bis}) has a symplectic structure, that is, that the restriction of the symplectic form $\sigma$ to $S$ is non-degenerate. We recall (see \cite{hps}) that this assumption is always fulfilled when the symbol $q$ fulfills (\ref{inf1}) and an assumption of partial ellipticity on its singular space $S$,
$$(x,\xi) \in S, \ q(x,\xi)=0 \Rightarrow (x,\xi)=0.$$
By denoting now $k_0$ the smallest integer $0 \leq k_0 \leq 2n-1$, such that
\begin{equation}\label{h1biskrist2bis}\inc
S=\Big(\bigcap_{j=0}^{k_0}\textrm{Ker}\big[\textrm{Re }F(\textrm{Im }F)^j \big]\Big) \cap \rr^{2n}, \num
\end{equation}
one can generalize Theorem~\ref{theorem1} as follows:

\bigskip

\begin{theorem}\label{theorem2}
Consider a quadratic operator $q(x,\xi)^w$ whose Weyl symbol  
\begin{eqnarray*}
q : \rr_x^n \times \rr_{\xi}^n &\rightarrow& \cc\\
 (x,\xi) & \mapsto & q(x,\xi),
\end{eqnarray*}
is a complex-valued quadratic form fulfilling \emph{(\ref{inf1})}. When its singular space $S$ has a symplectic structure then the operator $q(x,\xi)^w$ is subelliptic in any direction of $S^{\sigma \perp}$ in the sense that, if $(x',\xi')$ are some linear symplectic coordinates on $S^{\sigma \perp}$ then we have 
\begin{equation}\label{kristen5}\inc
\exists C>0, \forall u \in D(q), \  \big\|\big(\langle(x',\xi')\rangle^{2/(2k_0+1)}\big)^w u\big\|_{L^2} \leq C\big(\|q(x,\xi)^w u\|_{L^2}+\|u\|_{L^2}\big), \num
\end{equation}
where  $k_0$ stands for the smallest integer $0 \leq k_0 \leq 2n-1$ such that \emph{(\ref{h1biskrist2bis})} is fulfilled, and $\langle(x',\xi')\rangle=(1+|x'|^2+|\xi'|^2)^{1/2}$.
\end{theorem}

\bigskip

As we will see in the following, Theorem~\ref{theorem2} will be deduced from a simple adaptation of the analysis led in the proof of Theorem~\ref{theorem1}.

\bigskip

\noindent
\textit{Acknowledgements.} The author is particularly grateful to M.~Hitrik and N.~Lerner for very enriching comments and remarks about this work.

\section{Proof of Theorem~\ref{theorem1}}
\init

In the following, we shall use the notation $S_{\Omega}\big(m(X)^r,m(X)^{-2s}dX^2\big)$, where $\Omega$ is an open set in $\rr^{2n}$, $r,s \in \rr$ and $m \in C^{\infty}(\Omega,\rr_+^*)$, to stand for 
the class of symbols $a$ verifying
$$a \in C^{\infty}(\Omega), \ \forall  \alpha \in \nn^{2n}, \exists C_{\alpha}>0, \ |\partial_X^{\alpha}a(X)| \leq C_{\alpha} m(X)^{r-s|\alpha|}, \ X \in \Omega.$$
In the case where $\Omega=\rr^{2n}$, we shall drop for simplicity the index $\Omega$ in the notation.
We shall also use the notations $f \lesssim g$ and $f \sim g$, on $\Omega$, for respectively the estimates $\exists C>0$, $f \leq Cg$ and, $f \lesssim g$ and $g \lesssim f$, on $\Omega$.

The proof of Theorem~\ref{theorem1} will rely on the following key proposition. 
Considering
\begin{eqnarray*}
q : \rr_x^n \times \rr_{\xi}^n &\rightarrow& \cc\\
 (x,\xi) & \mapsto & q(x,\xi),
\end{eqnarray*}
a complex-valued quadratic form with a non-negative real part
\begin{equation}\label{giu1}\inc
\textrm{Re }q(x,\xi) \geq 0, \ (x,\xi) \in \rr^{2n}, n \in \nn^*, \num
\end{equation}
we assume that there exist a positive integer $m \in \nn^*$ and an open set $\Omega_0$ in $\rr^{2n}$ such that the following sum of $m+1$ non-negative quadratic forms satisfies
\begin{equation}\label{giu2}\inc  
\exists c_0>0, \forall X \in \Omega_0, \ \sum_{j=0}^{m}{\textrm{Re }q\big((\textrm{Im }F)^j X\big)} \geq c_0 |X|^2, \num
\end{equation}
where the notation $\textrm{Im }F$ stands for the imaginary part of the Hamilton map $F$ associated to the quadratic form $q$, then one can build a bounded weight function with the following properties:

\bigskip

\begin{proposition}\label{prop1}
If $q$ is a complex-valued quadratic form on $\rr^{2n}$ verifying \emph{(\ref{giu1})} and \emph{(\ref{giu2})} then there exists a real-valued weight function 
$$g \in S_{\Omega_0}\big(1,\langle X \rangle^{-\frac{2}{2m+1}}dX^2\big),$$ such that 
\begin{equation}\label{giu3}\inc
\exists c_1, c_2>0, \forall X \in \Omega_0, \ \emph{\textrm{Re }}q(X)+c_1H_{\emph{\textrm{Im}}q\ } g(X) +1 \geq c_2 \langle X \rangle^{\frac{2}{2m+1}}, \num
\end{equation}
where the notation $H_{\emph{\textrm{Im}}q}$ stands for the Hamilton vector field of the imaginary part of~$q$.
\end{proposition}

\bigskip

The construction of this weight function will really be the core of this paper. Its proof, which is technical, is given in Section~\ref{proofprop1}. Let us mention that because of its  
simple properties, this weight function may also be of further interest for future studies of doubly characteristic pseudodifferential operators with principal symbols whose Hessians at critical points fulfill (\ref{inf1}) and (\ref{h1biskrist1}).

Before proving this proposition, we shall explain how we can deduce Theorem~\ref{theorem1} from it. In doing so, we shall use as previously mentioned a multiplier method inspired from the work \cite{HeSjSt} of F.~H\'erau, J.~Sj\"ostrand and C.~Stolk about Fokker-Planck operators. In their analysis, they are led to establish a similar estimate as (\ref{kristen1}) in the case where the non-negative integer $k_0$ in Theorem~\ref{theorem1} is equal to 1. One can indeed check that their subelliptic assumption for their symbols at critical points, say here $X_0=0$,
$$\exists \eps_0>0, \ \textrm{Re }p(X)+\eps_0H_{\textrm{Im}p}^2\textrm{Re }p(X) \sim |X|^2,$$
is equivalent to the fact that their Hessians in these points fulfill (\ref{giu2}) with $m=1$ and $\Omega_0=\rr^{2n}$. In order to define our multiplier, we shall use the Wick quantization of the weight function given by Proposition~\ref{prop1}. The definition of the Wick quantization and some elements of Wick calculus, we need here, are recalled in the appendix (Section~\ref{wick}).

\medskip

To check that we can actually deduce Theorem~\ref{theorem1} from Proposition~\ref{prop1}, we begin by
considering a complex-valued quadratic form $q$ on $\rr^{2n}$, $n \geq 1$, with a non-negative real part and a zero singular space
$S=\{0\}.$ We know from (\ref{h1biskrist2}) that one can find a smallest integer $0 \leq k_0 \leq 2n-1$ such that
\begin{equation}\label{giu5}\inc
\Big(\bigcap_{j=0}^{k_0}\textrm{Ker}\big[\textrm{Re }F(\textrm{Im }F)^j \big]\Big) \cap \rr^{2n}=\{0\}.\num
\end{equation}
We then notice, as in~\cite{hps} and~\cite{mz}, that (\ref{giu5}) induces that the following sum of $k_0+1$ non-negative quadratic forms 
\begin{equation}\label{giu6}\inc  
\exists c_0>0, \forall X \in \rr^{2n}, \ r(X)=\sum_{j=0}^{k_0}{\textrm{Re }q\big((\textrm{Im }F)^j X\big)} \geq c_0 |X|^2, \num
\end{equation}
is a positive definite quadratic form. Let us indeed consider $X_0 \in \rr^{2n}$ such that $r(X_0)=0$. 
Then, the non-negativity of the quadratic form $\textrm{Re } q$ induces that for all $j=0,...,k_0$,
\begin{equation}\label{giu7}\inc
\textrm{Re }q\big((\textrm{Im } F)^j X_0\big)=0. \num
\end{equation}
By denoting by $\textrm{Re }q(X;Y)$ the polar form associated to $\textrm{Re }q$, we deduce from the Cauchy-Schwarz inequality, (\ref{10}) and (\ref{giu7}) that
for all $j=0,...,k_0$ and $Y \in \rr^{2n}$, 
\begin{align*}
\big|\textrm{Re }q\big(Y;(\textrm{Im }F)^j X_0\big)\big|^2= & \ \big|\sigma\big(Y,\textrm{Re }F (\textrm{Im }F)^j X_0\big)\big|^2 \\
 \leq & \ \textrm{Re }q(Y) \ \textrm{Re }q\big((\textrm{Im } F)^j X_0\big) =0. 
\end{align*} 
It follows that for all $j=0,...,k_0$ and $Y \in \rr^{2n}$, 
$$\sigma\big(Y,\textrm{Re }F (\textrm{Im }F)^j X_0\big)=0,$$
which implies that for all $j=0,...,k_0$, 
\begin{equation}\label{2.3.100}\inc
\textrm{Re }F(\textrm{Im }F)^j X_0=0, \num 
\end{equation}
since $\sigma$ is non-degenerate. We finally deduce (\ref{giu6}) from (\ref{giu5}).

In the case where $k_0=0$, the quadratic form $\textrm{Re }q$ is positive definite. This implies that the quadratic form $q$ is elliptic. As previously mentioned, the result of Theorem~\ref{theorem1} is in this case a straightforward consequence of classical results about elliptic quadratic differential operators recalled in (\ref{lay30}).

We can therefore assume in the following that $k_0 \geq 1$ and find from Proposition~\ref{prop1} a real-valued weight function 
\begin{equation}\label{giu8}\inc
g \in S\big(1,\langle X \rangle^{-\frac{2}{2k_0+1}}dX^2\big), \num
\end{equation}
such that 
\begin{equation}\label{giu9}\inc
\exists c_1, c_2>0, \forall X \in \rr^{2n}, \ \textrm{Re }q(X)+c_1H_{\textrm{Im}q\ } g(X)+1 \geq c_2 \langle X \rangle^{\frac{2}{2k_0+1}}. \num
\end{equation}
For $0<\eps \leq 1$, we consider the multiplier defined in the Wick quantization by the symbol $1-\eps g$. We recall that the definition of the Wick quantization and some elements of Wick calculus are recalled in Section~\ref{wick}. It follows from (\ref{giu8}), (\ref{lay0}), (\ref{lay1}), (\ref{lay2}) and the Cauchy-Schwarz inequality that 
\inc
\begin{multline*}\label{giu18}
\textrm{Re}\big(q^{\textrm{Wick}} u, (1-\eps g)^{\textrm{Wick}}u\big)=  \big(\textrm{Re}\big((1-\eps g)^{\textrm{Wick}} q^{\textrm{Wick}}\big) u,u\big) \num \\
\leq  \|1-\eps g\|_{L^{\infty}(\rr^{2n})}\|q^{\textrm{Wick}}u\|_{L^2}\|u\|_{L^2}  
\lesssim  \|q^{\textrm{Wick}}u\|_{L^2}^2 +\|u\|_{L^2}^2 \lesssim  \|\tilde{q}^{w}u\|_{L^2}^2 +\|u\|_{L^2}^2,
\end{multline*}
where  
\begin{equation}\label{lay5}\inc
\tilde{q}(x,\xi)=q\Big(x,\frac{\xi}{2\pi}\Big),\num
\end{equation}
because the operator $(1-\eps g)^{\textrm{Wick}}$ whose Wick symbol is real-valued, is formally selfadjoint. Indeed, the symbol $r(q)$ defined in (\ref{lay2}) is here just constant since $q$ is a quadratic form. The factor $2\pi$ in (\ref{lay5}) comes from the difference of normalizations chosen between (\ref{3}) and (\ref{lay3}) (See remark in Section~\ref{wick}). Since from (\ref{lay4}),
$$(1-\eps g)^{\textrm{Wick}} q^{\textrm{Wick}} =\Big{[}(1-\eps g)q+\frac{\eps}{4 \pi} \nabla g. \nabla q-\frac{\eps}{4i \pi}\{g,q\} \Big{]}^{\textrm{Wick}}+S,$$
with $\|S\|_{\mathcal{L}(L^2(\rr^n))} \lesssim 1$, we obtain from the fact real Hamiltonians get quantized in the Wick quantization by formally selfadjoint operators that 
$$\textrm{Re}\big((1-\eps g)^{\textrm{Wick}} q^{\textrm{Wick}}\big)=\Big{[}(1-\eps g)\textrm{Re }q+\frac{\eps}{4 \pi} \nabla g. \nabla \textrm{Re }q+\frac{\eps}{4 \pi}H_{\textrm{Im}q}\ g \Big{]}^{\textrm{Wick}}
+\textrm{Re }S,$$
because $g$ is a real-valued symbol. Since $\textrm{Re }q \geq 0$ and $g \in L^{\infty}(\rr^n)$, we can choose the positive parameter $\eps$ sufficiently small such that 
$$\forall X \in \rr^{2n}, \ 1-\eps g(X) \geq \frac{1}{2},$$
in order to deduce from (\ref{giu9}), (\ref{giu18}) and (\ref{lay0.5}) that 
\begin{equation}\label{lay11}\inc
\big((\langle X \rangle^{\frac{2}{2k_0+1}})^{\textrm{Wick}}u,u\big) \lesssim  \|\tilde{q}^{w}u\|_{L^2}^2 +\|u\|_{L^2}^2+ 
\big|\big((\nabla g. \nabla \textrm{Re }q)^{\textrm{Wick}}u,u\big)\big|, \num
\end{equation}
because from (\ref{lay0.1}) and (\ref{lay0.2}), $1^{\textrm{Wick}}=\textrm{Id}.$

By denoting $\tilde{X}=\big(x,\xi/(2\pi)\big)$ and $\textrm{Op}^w\big(S(1,dX^2)\big)$ the operators obtained by the Weyl quantization of symbols in the class $S(1,dX^2)$, it follows from (\ref{lay1}), (\ref{lay2}) and usual results of symbolic calculus that 
\begin{equation}\label{lay43}\inc
\big(\langle X \rangle^{\frac{2}{2k_0+1}}\big)^{\textrm{Wick}} - \big(\langle \tilde{X} \rangle^{\frac{2}{2k_0+1}}\big)^{w} \in \textrm{Op}^w\big(S(1,dX^2)\big) \num
\end{equation}
and
\begin{equation}\label{lay44}\inc
\big(\langle \tilde{X} \rangle^{\frac{1}{2k_0+1}}\big)^{w}\big(\langle \tilde{X} \rangle^{\frac{1}{2k_0+1}}\big)^{w}- \big(\langle \tilde{X} \rangle^{\frac{2}{2k_0+1}}\big)^{w} \in \textrm{Op}^w\big(S(1,dX^2)\big),\num
\end{equation}
since $k_0 \geq 0$. By using that 
$$\big(\big(\langle \tilde{X} \rangle^{\frac{1}{2k_0+1}}\big)^{w}\big(\langle \tilde{X} \rangle^{\frac{1}{2k_0+1}}\big)^{w}u,u\big)=\big\|\big(\langle \tilde{X} \rangle^{\frac{1}{2k_0+1}}\big)^{w}u\big\|_{L^2}^2,$$
we therefore deduce from (\ref{lay11}) and the Calder\'on-Vaillancourt theorem that 
\begin{equation}\label{lay12}\inc
\big\|\big(\langle \tilde{X} \rangle^{\frac{1}{2k_0+1}}\big)^{w}u\big\|_{L^2}^2 \lesssim  \|\tilde{q}^{w}u\|_{L^2}^2 +\|u\|_{L^2}^2+ 
\big|\big((\nabla g. \nabla \textrm{Re }q)^{\textrm{Wick}}u,u\big)\big|. \num
\end{equation}
Then, we get from (\ref{giu8}) and (\ref{lay0.5}) that 
\begin{equation}\label{lay13}\inc
\big|\big((\nabla g. \nabla \textrm{Re }q)^{\textrm{Wick}}u,u\big)\big| \lesssim \big(|\nabla \textrm{Re }q|^{\textrm{Wick}}u,u\big). \num
\end{equation}
Recalling now the well-known inequality 
\begin{equation}\label{giu00.1}\inc
|f'(x)|^2 \leq 2 f(x) \|f''\|_{L^{\infty}(\rr)}, \num
\end{equation}
fulfilled by any non-negative smooth function with bounded second derivative, we deduce from another use of (\ref{lay0.5}) that 
\begin{equation}\label{lay20}\inc
\big(|\nabla \textrm{Re }q|^{\textrm{Wick}}u,u\big) \lesssim \big(((\textrm{Re }q)^{\frac{1}{2}})^{\textrm{Wick}}u,u\big) \lesssim \big((1+ \textrm{Re }q)^{\textrm{Wick}}u,u\big), \num
\end{equation}
since $\textrm{Re }q$ is a non-negative quadratic form and that 
$$2(\textrm{Re }q)^{\frac{1}{2}} \leq 1+\textrm{Re }q.$$ 
By using the same arguments as in (\ref{giu18}), we obtain that 
\begin{multline*}
\big((1+ \textrm{Re }q)^{\textrm{Wick}}u,u\big)=\big((\textrm{Re }q)^{\textrm{Wick}}u,u\big)+\|u\|_{L^2}^2=\textrm{Re}(q^{\textrm{Wick}}u,u)+\|u\|_{L^2}^2\\
\leq \|q^{\textrm{Wick}}u\|_{L^2}\|u\|_{L^2}+\|u\|_{L^2}^2
\lesssim  \|q^{\textrm{Wick}}u\|_{L^2}^2 +\|u\|_{L^2}^2 \lesssim  \|\tilde{q}^{w}u\|_{L^2}^2 +\|u\|_{L^2}^2.
\end{multline*}
It therefore follows from (\ref{lay12}), (\ref{lay13}) and (\ref{lay20}) that 
\begin{equation}\label{lay21}\inc
\big\|\big(\langle \tilde{X} \rangle^{\frac{1}{2k_0+1}}\big)^{w}u\big\|_{L^2}^2 \lesssim  \|\tilde{q}^{w}u\|_{L^2}^2 +\|u\|_{L^2}^2. \num
\end{equation}
In order to improve the estimate (\ref{lay21}), we carefully resume our previous analysis and notice that our previous reasoning has in fact established that
\begin{align*}
\big\|\big(\langle \tilde{X} \rangle^{\frac{1}{2k_0+1}}\big)^{w}u\big\|_{L^2}^2 \lesssim & \  \big|\textrm{Re}\big(q^{\textrm{Wick}} u, (1-\eps g)^{\textrm{Wick}}u\big)\big|+ 
\big|\big((\nabla g. \nabla \textrm{Re }q)^{\textrm{Wick}}u,u\big)\big|+\|u\|_{L^2}^2 \\
\lesssim & \  \big|\textrm{Re}\big(q^{\textrm{Wick}} u, (1-\eps g)^{\textrm{Wick}}u\big)\big|+ 
|\textrm{Re}(q^{\textrm{Wick}} u, u)|+\|u\|_{L^2}^2\\
\lesssim & \  \big|\textrm{Re}\big(\tilde{q}^{w} u, (1-\eps g)^{\textrm{Wick}}u\big)\big|+ 
|\textrm{Re}(\tilde{q}^{w} u, u)|+\|u\|_{L^2}^2, 
\end{align*}
because $(1-\eps g)^{\textrm{Wick}}$ is a bounded operator on $L^2(\rr^n)$,
\begin{equation}\label{linn1}\inc
\|(1-\eps g)^{\textrm{Wick}}\|_{\mathcal{L}(L^2)} \leq \|1-\eps g\|_{L^{\infty}(\rr^{2n})}.\num
\end{equation}
By applying this estimate to $\big(\langle \tilde{X} \rangle^{\frac{1}{2k_0+1}}\big)^{w}u$, we deduce from (\ref{lay44}) and the Calder\'on-Vaillancourt theorem that 
\inc\begin{multline*}\label{lay45}
\big\|\big(\langle \tilde{X} \rangle^{\frac{2}{2k_0+1}}\big)^{w}u\big\|_{L^2}^2 \lesssim  \Big|\textrm{Re}\Big(\tilde{q}^{w} \big(\langle \tilde{X} \rangle^{\frac{1}{2k_0+1}}\big)^{w} u, 
(1-\eps g)^{\textrm{Wick}}\big(\langle \tilde{X} \rangle^{\frac{1}{2k_0+1}}\big)^{w}u\Big)\Big|\\ + 
\Big|\textrm{Re}\Big(\tilde{q}^{w} \big(\langle \tilde{X} \rangle^{\frac{1}{2k_0+1}}\big)^{w}u, \big(\langle \tilde{X} \rangle^{\frac{1}{2k_0+1}}\big)^{w}u\Big)\Big|+\big\|\big(\langle \tilde{X} \rangle^{\frac{1}{2k_0+1}}\big)^{w}u\big\|_{L^2}^2+\|u\|_{L^2}^2. \num
\end{multline*}
Then, by noticing that the commutator 
\begin{equation}\label{lay50}\inc
\big[\tilde{q}^w,\big(\langle \tilde{X} \rangle^{\frac{1}{2k_0+1}}\big)^{w}\big] \in \textrm{Op}^w\big(S\big(\langle X \rangle^{\frac{1}{2k_0+1}},\langle X \rangle^{-2}dX^2\big)\big), \num
\end{equation}
because $\tilde{q}$ is a quadratic form, and that 
\begin{equation}\label{lay51}\inc
\big(\langle \tilde{X} \rangle^{-\frac{1}{2k_0+1}}\big)^{w}\big(\langle \tilde{X} \rangle^{\frac{1}{2k_0+1}}\big)^{w}- \textrm{Id} \in \textrm{Op}^w\big(S(\langle X \rangle^{-2},\langle X \rangle^{-2}dX^2)\big),\num
\end{equation}
we deduce from standard results of symbolic calculus and the Calder\'on-Vaillancourt theorem that 
\inc\begin{align*}\label{lay52}
\big\|\big[\tilde{q}^w,\big(\langle \tilde{X} \rangle^{\frac{1}{2k_0+1}}\big)^{w}\big]u\big\|_{L^2} \lesssim & \ \big\|\big[\tilde{q}^w,\big(\langle \tilde{X} \rangle^{\frac{1}{2k_0+1}}\big)^{w}\big]\big(\langle \tilde{X} \rangle^{-\frac{1}{2k_0+1}}\big)^{w}\big(\langle \tilde{X} \rangle^{\frac{1}{2k_0+1}}\big)^{w}u\big\|_{L^2}+ \|u\|_{L^2}\\
\lesssim & \ \big\|\big(\langle \tilde{X} \rangle^{\frac{1}{2k_0+1}}\big)^{w}u\big\|_{L^2}+\|u\|_{L^2}. \num
\end{align*}
By introducing this commutator, we get from
the Cauchy-Schwarz inequality and (\ref{lay52}) that 
\begin{multline*}
\Big|\textrm{Re}\Big(\tilde{q}^{w} \big(\langle \tilde{X} \rangle^{\frac{1}{2k_0+1}}\big)^{w}u, \big(\langle \tilde{X} \rangle^{\frac{1}{2k_0+1}}\big)^{w}u\Big)\Big| \lesssim 
\Big|\textrm{Re}\Big(\tilde{q}^{w} u, \big(\langle \tilde{X} \rangle^{\frac{1}{2k_0+1}}\big)^{w}\big(\langle \tilde{X} \rangle^{\frac{1}{2k_0+1}}\big)^{w}u\Big)\Big|
\\
+\big\|\big(\langle \tilde{X} \rangle^{\frac{1}{2k_0+1}}\big)^{w}u\big\|_{L^2}^2+\|u\|_{L^2}^2.
\end{multline*}
By using that another use of the Cauchy-Schwarz inequality and the Calder\'on-Vaillancourt theorem with (\ref{lay44}) gives that 
$$\Big|\textrm{Re}\Big(\tilde{q}^{w} u, \big(\langle \tilde{X} \rangle^{\frac{1}{2k_0+1}}\big)^{w}\big(\langle \tilde{X} \rangle^{\frac{1}{2k_0+1}}\big)^{w}u\Big)\Big| \lesssim
\|\tilde{q}^{w} u\|_{L^2}\big\|\big(\langle \tilde{X} \rangle^{\frac{2}{2k_0+1}}\big)^{w}u\big\|_{L^2}+\|\tilde{q}^{w} u\|_{L^2}\|u\|_{L^2},$$
we deduce from (\ref{lay21}) and the previous estimate that
\begin{align*}
& \ \Big|\textrm{Re}\Big(\tilde{q}^{w} \big(\langle \tilde{X} \rangle^{\frac{1}{2k_0+1}}\big)^{w}u, \big(\langle \tilde{X} \rangle^{\frac{1}{2k_0+1}}\big)^{w}u\Big)\Big| \\ 
\lesssim & \
\|\tilde{q}^{w} u\|_{L^2}\big\|\big(\langle \tilde{X} \rangle^{\frac{2}{2k_0+1}}\big)^{w}u\big\|_{L^2}
+\|\tilde{q}^{w} u\|_{L^2}^2+\|u\|_{L^2}^2.
\end{align*}
By using again the Cauchy-Schwarz inequality, (\ref{lay21}), (\ref{linn1}), (\ref{lay45}) and (\ref{lay52}), this estimate implies that 
\inc\begin{align*}\label{lay46}
& \ \big\|\big(\langle \tilde{X} \rangle^{\frac{2}{2k_0+1}}\big)^{w}u\big\|_{L^2}^2 \lesssim  \Big|\textrm{Re}\Big(\big[\tilde{q}^{w}, \big(\langle \tilde{X} \rangle^{\frac{1}{2k_0+1}}\big)^{w}\big] u, 
(1-\eps g)^{\textrm{Wick}}\big(\langle \tilde{X} \rangle^{\frac{1}{2k_0+1}}\big)^{w}u\Big)\Big| \num \\ 
+& \ \Big|\textrm{Re}\Big(\tilde{q}^{w}  u, 
\big(\langle \tilde{X} \rangle^{\frac{1}{2k_0+1}}\big)^{w}(1-\eps g)^{\textrm{Wick}}\big(\langle \tilde{X} \rangle^{\frac{1}{2k_0+1}}\big)^{w}u\Big)\Big|
+\|\tilde{q}^{w} u\|_{L^2}^2+\|u\|_{L^2}^2 \\
\lesssim & \ \Big|\textrm{Re}\Big(\tilde{q}^{w}  u, 
\big(\langle \tilde{X} \rangle^{\frac{1}{2k_0+1}}\big)^{w}(1-\eps g)^{\textrm{Wick}}\big(\langle \tilde{X} \rangle^{\frac{1}{2k_0+1}}\big)^{w}u\Big)\Big|
+\|\tilde{q}^{w} u\|_{L^2}^2+\|u\|_{L^2}^2\\
\lesssim & \ \|\tilde{q}^{w}  u\|_{L^2}\big\| 
\big(\langle \tilde{X} \rangle^{\frac{1}{2k_0+1}}\big)^{w}(1-\eps g)^{\textrm{Wick}}\big(\langle \tilde{X} \rangle^{\frac{1}{2k_0+1}}\big)^{w}u\big\|_{L^2}
+\|\tilde{q}^{w} u\|_{L^2}^2+\|u\|_{L^2}^2,
\end{align*}
because we get from (\ref{linn1}) and (\ref{lay52}) that 
\begin{multline*}
\Big|\textrm{Re}\Big(\big[\tilde{q}^{w}, \big(\langle \tilde{X} \rangle^{\frac{1}{2k_0+1}}\big)^{w}\big] u, 
(1-\eps g)^{\textrm{Wick}}\big(\langle \tilde{X} \rangle^{\frac{1}{2k_0+1}}\big)^{w}u\Big)\Big| \lesssim 
\big\|\big(\langle \tilde{X} \rangle^{\frac{1}{2k_0+1}}\big)^{w}u\big\|_{L^2}^2\\ +\big\|\big(\langle \tilde{X} \rangle^{\frac{1}{2k_0+1}}\big)^{w}u\big\|_{L^2}\|u\|_{L^2}.
\end{multline*}
Notice now that (\ref{giu8}), (\ref{lay1bis}) and (\ref{lay2bis}) imply that 
$$\big[\big(\langle \tilde{X} \rangle^{\frac{1}{2k_0+1}}\big)^{w},(1-\eps g)^{\textrm{Wick}}\big] \in \textrm{Op}^w\big(S(1,dX^2)\big),$$
since $(1-\eps g)^{\textrm{Wick}}=\tilde{g}^w$, with $\tilde{g} \in S(1,dX^2)$ and $k_0 \geq 0$. 
By introducing this new commutator, we deduce from the Calder\'on-Vaillancourt theorem, (\ref{lay44}), (\ref{lay21}) and (\ref{linn1}) that
\begin{align*}
& \ \big\| \big(\langle \tilde{X} \rangle^{\frac{1}{2k_0+1}}\big)^{w}(1-\eps g)^{\textrm{Wick}}\big(\langle \tilde{X} \rangle^{\frac{1}{2k_0+1}}\big)^{w}u\big\|_{L^2} \\
\lesssim & \ \big\|\big(\langle \tilde{X} \rangle^{\frac{1}{2k_0+1}}\big)^{w}u\big\|_{L^2}
+\big\| (1-\eps g)^{\textrm{Wick}}\big(\langle \tilde{X} \rangle^{\frac{1}{2k_0+1}}\big)^{w}\big(\langle \tilde{X} \rangle^{\frac{1}{2k_0+1}}\big)^{w}u\big\|_{L^2} \\
\lesssim & \ \big\|\big(\langle \tilde{X} \rangle^{\frac{1}{2k_0+1}}\big)^{w}u\big\|_{L^2}
+\big\|\big(\langle \tilde{X} \rangle^{\frac{1}{2k_0+1}}\big)^{w}\big(\langle \tilde{X} \rangle^{\frac{1}{2k_0+1}}\big)^{w}u\big\|_{L^2} \\\lesssim & \ \big\|\big(\langle \tilde{X} \rangle^{\frac{2}{2k_0+1}}\big)^{w}u\big\|_{L^2}+\big\|\big(\langle \tilde{X} \rangle^{\frac{1}{2k_0+1}}\big)^{w}u\big\|_{L^2}+\|u\|_{L^2}\\
\lesssim & \ \big\|\big(\langle \tilde{X} \rangle^{\frac{2}{2k_0+1}}\big)^{w}u\big\|_{L^2}+\|\tilde{q}^{w}u\|_{L^2}+\|u\|_{L^2}.
\end{align*}
Recalling (\ref{lay46}), we can then use this last estimate to obtain that 
\begin{equation}\label{lay60}\inc
\big\|\big(\langle \tilde{X} \rangle^{\frac{2}{2k_0+1}}\big)^{w}u\big\|_{L^2}^2 \lesssim \|\tilde{q}^{w}u\|_{L^2}^2+\|u\|_{L^2}^2. \num
\end{equation}
By finally noticing from the homogeneity of degree 2 of $\tilde{q}$ that we have 
$$(\tilde{q} \circ T)(x,\xi)=\frac{1}{2\pi}q(x,\xi),$$
if $T$ stands for the real linear symplectic transformation 
$$T(x,\xi)=\big((2\pi)^{-\frac{1}{2}}x,(2\pi)^{\frac{1}{2}}\xi\big),$$
we deduce from the symplectic invariance of the Weyl quantization (Theorem~18.5.9 in~\cite{hormander}) that  
$$\big\|\big(\langle X \rangle^{\frac{2}{2k_0+1}}\big)^{w}u\big\|_{L^2}^2 \lesssim  \|q^{w}u\|_{L^2}^2 +\|u\|_{L^2}^2,$$
which proves Theorem~\ref{theorem1}.

\section{Proof of Theorem~\ref{theorem2}}
\init

This section is devoted to the proof of Theorem~\ref{theorem2}. We begin by recalling that the symplectic invariance property of the Weyl quantization (Theorem~18.5.9 in \cite{hormander}) allows us to freely choose the linear symplectic coordinates $(x,\xi)$ in which we want to express our symbol $q$ in our proof of Theorem~\ref{theorem2}.  
Considering
\begin{eqnarray*}
q : \rr_x^n \times \rr_{\xi}^n &\rightarrow& \cc\\
 (x,\xi) & \mapsto & q(x,\xi),
\end{eqnarray*}
a complex-valued quadratic form with a non-negative real part
$$\textrm{Re }q(x,\xi) \geq 0, \ (x,\xi) \in \rr^{2n}, n \in \nn^*,$$
and assuming that its singular space $S$ has a symplectic structure, we deduce from Proposition~2.0.1 in~\cite{hps} that one can find some linear symplectic coordinates in $\rr^{2n}$,
$$(x,\xi)=(x',x'';\xi',\xi'') \in \rr^{2n},$$ 
with $(x',\xi')$ and $(x'',\xi'')$ some linear symplectic coordinates respectively in $S^{\sigma \perp}$ and $S$; such that we can write the symbol $q$ as the sum of two quadratic forms 
\begin{equation}\label{emi4}\inc
q(x,\xi)=q_1(x',\xi')+iq_2(x'',\xi''),\num
\end{equation}
where $q_1$ is a complex-valued quadratic form on $\rr^{2n'}$ with a non-negative real part and $q_2$ is a real-valued quadratic form on $\rr^{2n''}$. More precisely, we proved in \cite{hps} (Proposition~2.0.1) that the spaces $S$ and $S^{\sigma \perp}$ are stable by the real and imaginary parts of the Hamilton map $F$ of the symbol $q$; and that the two quadratic forms $q_1$ and $q_2$ are actually equal to
$$q_1(x',\xi')=\sigma\big((x',\xi'),F|_{S^{\sigma \perp}}(x',\xi')\big) \textrm{ and } iq_2(x'',\xi'')=\sigma\big((x'',\xi''),F|_{S}(x'',\xi'')\big).$$
By denoting $F_1=F|_{S^{\sigma \perp}}$ the Hamilton map of $q_1$,  we first check that (\ref{h1biskrist2bis}) implies that the non-negative quadratic form
\begin{equation}\label{emi1}\inc
r(X')=\sum_{j=0}^{k_0}\textrm{Re }q_1\big((\textrm{Im }F_1)^jX'\big), \num
\end{equation}
is actually positive definite on $S^{\sigma \perp}$. Indeed, consider $X_0' \in S^{\sigma \perp}$ such that $r(X_0')=0$. As in (\ref{2.3.100}), it follows that 
$\textrm{Re }F_1(\textrm{Im }F_1)^jX_0'=0$ for all $0 \leq j \leq k_0$, which according to (\ref{h1biskrist2bis}), implies that 
$X_0' \in S \cap S^{\sigma \perp} =\{0\}$.

Let us first consider the case where $k_0 \geq 1$. As in the proof of Theorem~\ref{theorem1}, one can find from (\ref{emi1}) and Proposition~\ref{prop1} a real-valued weight function in the variables $X'=(x',\xi') \in S^{\sigma \perp}$,
\begin{equation}\label{emi2}\inc
g_1 \in S\big(1,\langle X' \rangle^{-\frac{2}{2k_0+1}}dX'^2\big), \num
\end{equation}
such that 
\begin{equation}\label{emi3}\inc
\exists c_1, c_2>0, \forall X' \in S^{\sigma \perp}, \ \textrm{Re }q_1(X')+c_1H_{\textrm{Im}q_1\ } g_1(X')+1 \geq c_2 \langle X' \rangle^{\frac{2}{2k_0+1}}. \num
\end{equation}
When $k_0=0$, it is sufficient to just take $g_1=0$ to fulfill (\ref{emi3}). Then, as previously in (\ref{giu18}), one can use the multiplier defined in the Wick quantization by the symbol $1-\eps g_1$, for $0<\eps \leq 1$; and consider the quantity
$$\textrm{Re}\big(q^{\textrm{Wick}} u, (1-\eps g_1)^{\textrm{Wick}}u\big).$$
By noticing from (\ref{emi4}) that we have this time 
$$\textrm{Re}\big((1-\eps g_1)^{\textrm{Wick}} q^{\textrm{Wick}}\big)=\Big{[}(1-\eps g_1)\textrm{Re }q_1+\frac{\eps}{4 \pi} \nabla g_1. \nabla \textrm{Re }q_1+\frac{\eps}{4 \pi}H_{\textrm{Im}q_1}\ g_1 \Big{]}^{\textrm{Wick}}
+\textrm{Re }S_1,$$
with $\|\textrm{Re }S_1\|_{\mathcal{L}(L^2)} \lesssim 1$, since 
$$H_{\textrm{Im}q} \ g_1=H_{\textrm{Im}q_1} \ g_1,$$ 
because of the variables tensorization. Next, one can exactly resume our analysis led in the proof of Theorem~\ref{theorem1} from (\ref{lay11}) in order to finish the proof of Theorem~\ref{theorem2}.

\section{Proof of Proposition~\ref{prop1}}\label{proofprop1}
\init

We prove the proposition~\ref{prop1} by induction on the positive integer $m \geq 1$ appearing in (\ref{giu2}). Let $m \geq 1$,  we shall assume that the proposition~\ref{prop1} is fulfilled for any 
open set $\Omega_0$ of $\rr^{2n}$, when the positive integer in (\ref{giu2}) is strictly smaller than $m$.

In the following, we denote by $\psi$, $\chi$ and $w$ some $C^{\infty}(\rr,[0,1])$ functions respectively satisfying  
\begin{equation}\label{giu13}\inc
\psi=1 \textrm{ on } [-1,1], \ \textrm{supp } \psi \subset [-2,2], \num
\end{equation} 
\begin{equation}\label{giu14}\inc
\chi=1 \textrm{ on } \{x \in \rr : 1 \leq |x| \leq 2 \}, \ \textrm{supp } \chi \subset \{x \in \rr :  1/2 \leq |x| \leq 3\}, \num
\end{equation}
and
\begin{equation}\label{giu15}\inc
w=1 \textrm{ on } \{x \in \rr : |x| \geq 2 \}, \ \textrm{supp } w \subset \{x \in \rr :  |x| \geq 1\}. \num
\end{equation}
More generically, we shall denote by $\psi_j$, $\chi_j$ and $w_j$, $j \in \nn$, some other $C^{\infty}(\rr,[0,1])$ functions satisfying similar properties as respectively $\psi$, $\chi$ and $w$ with possibly different choices for the positive numerical values which define their support localizations.

Let $\Omega_0$ be an open set of $\rr^{2n}$ such that (\ref{giu2}) is fulfilled. Considering the quadratic form 
\begin{equation}\label{giu11}\inc
r_k(X)=\textrm{Re }q\big((\textrm{Im }F)^{k-1}X;(\textrm{Im }F)^{k}X\big), \ k \in \nn^*, \num
\end{equation}
and defining
\begin{equation}\label{giu23}\inc
g_m(X)=\psi\Big(\textrm{Re }q\big((\textrm{Im }F)^{m-1}X\big)\langle X \rangle^{-\frac{2(2m-1)}{2m+1}}\Big)\langle X \rangle^{-\frac{4m}{2m+1}}r_m(X), \num
\end{equation}
where $\psi$ is the function defined in (\ref{giu13}), we get from Lemma~\ref{lem2} that 
\begin{align*}\label{giu24}\inc
& \ H_{\textrm{Im}q}\ g_m(X)=2\psi\Big(\textrm{Re }q\big((\textrm{Im }F)^{m-1}X\big)\langle X \rangle^{-\frac{2(2m-1)}{2m+1}}\Big)\frac{\textrm{Re }q\big((\textrm{Im }F)^mX\big)}{\langle X \rangle^{\frac{4m}{2m+1}}} \num \\
+ & \ 2\psi\Big(\textrm{Re }q\big((\textrm{Im }F)^{m-1}X\big)\langle X \rangle^{-\frac{2(2m-1)}{2m+1}}\Big)\frac{\textrm{Re }q\big((\textrm{Im }F)^{m-1}X;(\textrm{Im }F)^{m+1}X\big)}{\langle X \rangle^{\frac{4m}{2m+1}}}\\
+ & \ H_{\textrm{Im}q}\Big(\psi\Big(\textrm{Re }q\big((\textrm{Im }F)^{m-1}X\big)\ \langle X \rangle^{-\frac{2(2m-1)}{2m+1}}\Big)\Big)\frac{r_m(X)}{\langle X \rangle^{\frac{4m}{2m+1}}}\\
+ & \ \psi\Big(\textrm{Re }q\big((\textrm{Im }F)^{m-1}X\big)\ \langle X \rangle^{-\frac{2(2m-1)}{2m+1}}\Big)H_{\textrm{Im}q}\Big(\langle X \rangle^{-\frac{4m}{2m+1}}\Big)r_m(X).
\end{align*}
We first check that 
\begin{equation}\label{giu25}\inc
g_m \in S\big(1,\langle X \rangle^{-\frac{2(2m-1)}{2m+1}}dX^2\big).\num
\end{equation}
In order to verify this, we notice from Lemma~\ref{lem2.21} that the two quadratic forms 
\begin{equation}\label{lay80}\inc
\textrm{Re }q\big((\textrm{Im }F)^{m-1}X;(\textrm{Im }F)^{m}X\big) \textrm{ and } \textrm{Re }q\big((\textrm{Im }F)^{m-1}X;(\textrm{Im }F)^{m+1}X\big),\num
\end{equation} 
belong to the symbol class 
\begin{equation}\label{lay81}\inc
S_{\Omega}\big(\langle X \rangle^{\frac{4m}{2m+1}}, \langle X \rangle^{-\frac{2(2m-1)}{2m+1}}dX^2\big),\num
\end{equation}
for any open set $\Omega$ where $\textrm{Re }q\big((\textrm{Im }F)^{m-1}X\big) \lesssim \langle X \rangle^{\frac{2(2m-1)}{2m+1}}$. To check this, we just use in addition to Lemma~\ref{lem2.21} the obvious estimates
$$\textrm{Re }q\big((\textrm{Im }F)^{m}X\big)^{\frac{1}{2}} \lesssim \langle X \rangle \textrm{ and } \textrm{Re }q\big((\textrm{Im }F)^{m+1}X\big)^{\frac{1}{2}} \lesssim \langle X \rangle.$$
Moreover, since 
\begin{equation}\label{oc1}\inc
\langle X \rangle^{-\frac{4m}{2m+1}} \in S\big(\langle X \rangle^{-\frac{4m}{2m+1}},\langle X \rangle^{-2}dX^2\big),\num
\end{equation}
we obtain (\ref{giu25}) from (\ref{giu13}), (\ref{giu11}), (\ref{giu23}), (\ref{lay80}), (\ref{lay81}) and Lemma~\ref{lem3}.

Denoting respectively $A_1$, $A_2$, $A_3$ and $A_4$ the four terms appearing in the right hand side of (\ref{giu24}), we first notice from (\ref{giu13}), (\ref{lay80}), (\ref{lay81}), (\ref{oc1}) and Lemma~\ref{lem3} that 
\begin{equation}\label{giu27}\inc
A_2 \in  S\big(1,\langle X \rangle^{-\frac{2(2m-1)}{2m+1}}dX^2\big).\num
\end{equation}
Next, by using that 
$$\textrm{Im }q \in S\big(\langle X \rangle^2,\langle X \rangle^{-2}dX^2\big),$$
since $\textrm{Im }q$ is a quadratic form,
we get from (\ref{giu13}), (\ref{giu11}), (\ref{lay80}), (\ref{lay81}), (\ref{oc1}) and Lemma~\ref{lem3} that
\begin{equation}\label{lay82}\inc
A_3 \in S\big(\langle X \rangle^{\frac{2}{2m+1}},\langle X \rangle^{-\frac{2(2m-1)}{2m+1}}dX^2\big),\num
\end{equation}
since
$$H_{\textrm{Im}q}\Big(\psi\Big(\textrm{Re }q\big((\textrm{Im }F)^{m-1}X\big)\ \langle X \rangle^{-\frac{2(2m-1)}{2m+1}}\Big)\Big)\in S\big(\langle X \rangle^{\frac{2}{2m+1}},\langle X \rangle^{-\frac{2(2m-1)}{2m+1}}dX^2\big).$$
By using now that 
$$H_{\textrm{Im}q}\Big(\langle X \rangle^{-\frac{4m}{2m+1}}\Big) \in S\big(\langle X \rangle^{-\frac{4m}{2m+1}},\langle X \rangle^{-2}dX^2\big),$$
we finally obtain from another use of (\ref{giu13}), (\ref{giu11}), (\ref{lay80}), (\ref{lay81}) and Lemma~\ref{lem3} that 
\begin{equation}\label{lay83}\inc
A_4 \in S\big(1,\langle X \rangle^{-\frac{2(2m-1)}{2m+1}}dX^2\big).\num
\end{equation}
Since the term $A_3$ is supported in 
$$\textrm{supp } \psi'\Big(\textrm{Re }q\big((\textrm{Im }F)^{m-1}X\big)\langle X \rangle^{-\frac{2(2m-1)}{2m+1}}\Big),$$ 
we deduce from (\ref{giu24}), (\ref{giu27}), (\ref{lay82}) and (\ref{lay83}) that there exists $\chi_0$ a $C^{\infty}(\rr,[0,1])$ function satisfying similar properties as in (\ref{giu14}), with possibly different positive numerical values for its support localization, such that, $\exists c_1,c_2>0$, $\forall X \in \rr^{2n}$,
\begin{align*}\label{giu28}\inc
 & \ H_{\textrm{Im}q}\ g_m(X)+c_1+
c_2\chi_0\Big(\textrm{Re }q\big((\textrm{Im }F)^{m-1}X\big)\langle X \rangle^{-\frac{2(2m-1)}{2m+1}}\Big)\langle X \rangle^{\frac{2}{2m+1}} \num \\
 & \ \geq 2 \psi\Big(\textrm{Re }q\big((\textrm{Im }F)^{m-1}X\big)\langle X \rangle^{-\frac{2(2m-1)}{2m+1}}\Big)\frac{\textrm{Re }q\big((\textrm{Im }F)^mX\big)}{\langle X \rangle^{\frac{4m}{2m+1}}}. 
\end{align*}
Recalling (\ref{giu2}), one can find some positive constants $c_3, c_4>0$ such that
\begin{equation}\label{giu29}\inc  
\sum_{j=0}^{m-1}{\textrm{Re }q\big((\textrm{Im }F)^j X\big)} \geq c_3 |X|^2, \num
\end{equation}on the open set 
\begin{equation}\label{oc2}\inc
\Omega_1=\big\{X \in \rr^{2n} : \textrm{Re }q\big((\textrm{Im }F)^{m}X\big) < c_4|X|^2\big\} \cap \Omega_0. \num
\end{equation}
When $m \geq 2$, one can find according to our induction hypothesis a real-valued function 
\begin{equation}\label{giu29.1}\inc
\tilde{g}_m \in S_{\Omega_1}\big(1,\langle X \rangle^{-\frac{2}{2m-1}}dX^2\big), \num
\end{equation}
such that 
\begin{equation}\label{giu30}\inc
\exists c_5>0, \forall X \in \Omega_1, \ \textrm{Re }q(X)+c_5H_{\textrm{Im}q\ } \tilde{g}_m(X) +1 \gtrsim  \langle X \rangle^{\frac{2}{2m-1}}. \num
\end{equation}
For convenience, we set in the following $\tilde{g}_1=0$ when $m=1$.  
By choosing suitably $\psi_0$ and $w_0$ some $C^{\infty}(\rr,[0,1])$ functions satisfying similar properties as the functions respectively defined in (\ref{giu13}) and (\ref{giu15}), with possibly different positive numerical values for their support localizations, such that 
\inc\begin{multline*}\label{bri2}
\textrm{supp } \psi_0\Big(\textrm{Re }q\big((\textrm{Im }F)^{m}X\big)|X|^{-2}\Big)w_0(X) \\ \subset \big\{X \in \rr^{2n} : \textrm{Re }q\big((\textrm{Im }F)^{m}X\big) < c_4|X|^2\big\},\num
\end{multline*}
and setting 
\begin{equation}\label{giu32}\inc
G_m(X)=g_m(X)+ \psi_0\Big(\textrm{Re }q\big((\textrm{Im }F)^{m}X\big) |X|^{-2}\Big)w_0(X)\tilde{g}_m(X), \ X \in \Omega_0, \num
\end{equation} 
we deduce from a straightforward adaptation of the Lemma~\ref{lem3} by recalling (\ref{giu13}) and (\ref{giu15}) that 
\begin{equation}\label{bri1}\inc
\psi_0\Big(\textrm{Re }q\big((\textrm{Im }F)^{m}X\big) |X|^{-2}\Big)w_0(X) \in S\big(1,\langle X \rangle^{-2}dX^2\big).\num
\end{equation}
According to (\ref{giu25}) and (\ref{giu29.1}), this implies that
\begin{equation}\label{giu33}\inc
G_1 \in S_{\Omega_0}\big(1,\langle X \rangle^{-\frac{2}{3}}dX^2\big) \textrm{ and } G_m \in S_{\Omega_0}\big(1,\langle X \rangle^{-\frac{2}{2m-1}}dX^2\big), \num
\end{equation}
when $m \geq 2$.
Since from (\ref{bri1}),
$$H_{\textrm{Im}q}\Big(\psi_0\Big(\textrm{Re }q\big((\textrm{Im }F)^{m}X\big)|X|^{-2}\Big)w_0(X) \Big) \in S\big(1,\langle X \rangle^{-2}dX^2\big),$$
because $\textrm{Im }q$ is a quadratic form, 
we first notice from (\ref{oc2}), (\ref{giu29.1}) and (\ref{bri2}) that 
$$H_{\textrm{Im}q}\Big(\psi_0\Big(\textrm{Re }q\big((\textrm{Im }F)^{m}X\big)|X|^{-2}\Big)w_0(X) \Big)\tilde{g}_m(X)
\in S_{\Omega_0}\big(1,\langle X \rangle^{-\frac{2}{2m-1}}dX^2\big),$$ 
and then deduce from (\ref{giu28}), (\ref{oc2}), (\ref{giu30}), (\ref{bri2}) and (\ref{giu32}) that there exist $c_6,c_7>0$ such that for all $X \in \Omega_0$,
\begin{align*}
& \ \textrm{Re }q(X)+c_6 H_{\textrm{Im}q}\ G_m(X)+1+
c_7\chi_0\Big(\textrm{Re }q\big((\textrm{Im }F)^{m-1}X\big)\langle X \rangle^{-\frac{2(2m-1)}{2m+1}}\Big)\langle X \rangle^{\frac{2}{2m+1}}  \\
& \ \gtrsim \psi\Big(\textrm{Re }q\big((\textrm{Im }F)^{m-1}X\big)\langle X \rangle^{-\frac{2(2m-1)}{2m+1}}\Big)\frac{\textrm{Re }q\big((\textrm{Im }F)^mX\big)}{\langle X \rangle^{\frac{4m}{2m+1}}} \\
& \ \qquad +\psi_0\Big(\textrm{Re }q\big((\textrm{Im }F)^{m}X\big)|X |^{-2}\Big)w_0(X)\langle X \rangle^{\frac{2}{2m-1}}, 
\end{align*}
when $m \geq 2$.
Since 
$$\langle X \rangle^{\frac{2}{2m-1}} \gtrsim \langle X \rangle^{\frac{2}{2m+1}} \textrm{ and }
\frac{\textrm{Re }q\big((\textrm{Im }F)^mX\big)}{\langle X \rangle^{\frac{4m}{2m+1}}} \gtrsim |X|^{\frac{2}{2m+1}},$$
when $\textrm{Re }q\big((\textrm{Im }F)^mX\big) \gtrsim  |X|^2$, 
we deduce from the previous estimate by distinguishing the regions in $\Omega_0$ where
$$\textrm{Re }q\big((\textrm{Im }F)^mX\big) \lesssim  |X|^2 \textrm{ and } \textrm{Re }q\big((\textrm{Im }F)^mX\big) \gtrsim  |X|^2,$$ 
according to the support of the function
$$\psi_0\Big(\textrm{Re }q\big((\textrm{Im }F)^{m}X\big)|X |^{-2}\Big),$$
that one can find  a $C^{\infty}(\rr,[0,1])$ function $w_1$ with the same kind of support as the function defined in (\ref{giu15}) such that 
\inc\begin{multline*}\label{giu35}
\exists c_8,c_9>0, \forall X \in \Omega_0, \  \textrm{Re }q(X) \num +c_8H_{\textrm{Im}q}\ G_m(X)\\ +c_9 w_1\Big(\textrm{Re }q\big((\textrm{Im }F)^{m-1}X\big)\langle X \rangle^{-\frac{2(2m-1)}{2m+1}}\Big)\langle X \rangle^{\frac{2}{2m+1}} 
+1 \gtrsim \langle X \rangle^{\frac{2}{2m+1}}, 
\end{multline*}
when $m \geq 2$. When $m=1$, 
we notice from (\ref{giu2}) that
\begin{equation}\label{bri3}\inc
\textrm{Re }q\big(\textrm{Im }FX\big) \gtrsim \langle X \rangle^{2},\num
\end{equation}
on any set where
\begin{equation}\label{bri4}\inc
|X|\geq c_{10} \textrm{ and } \textrm{Re }q(X) \leq  \langle X \rangle^{\frac{2}{3}},\num
\end{equation}
if the positive constant $c_{10}$ is chosen sufficiently large.  Moreover, since in this case $G_1=g_1$ and that $\textrm{Re }q \geq 0$, one can deduce from (\ref{giu13}), (\ref{giu15}), (\ref{giu28}), (\ref{bri3}) and (\ref{bri4}),
by distinguishing the regions in $\Omega_0$ where 
$$\textrm{Re }q(X) \lesssim \langle X \rangle^{\frac{2}{3}} \textrm{ and } \textrm{Re }q(X) \gtrsim \langle X \rangle^{\frac{2}{3}},$$
according to the support of the function
$$\psi\big(\textrm{Re }q(X)\langle X \rangle^{-\frac{2}{3}}\big),$$
that the estimate (\ref{giu35}) is also fulfilled in the case $m=1$. Continuing our study of the case where $m=1$, we notice from (\ref{giu15}) and $\textrm{Re }q \geq 0$, that one can estimate
$$w_1\big(\textrm{Re }q(X)\langle X \rangle^{-\frac{2}{3}}\big)\langle X \rangle^{\frac{2}{3}} \lesssim \textrm{Re }q(X),$$
for all $X \in \rr^{2n}$. It therefore follows that one can find $c_{11}>0$ such that for all $X \in \Omega_0$,
$$\textrm{Re }q(X) +c_{11}H_{\textrm{Im}q}\ G_1(X) +1 \gtrsim \langle X \rangle^{\frac{2}{3}},$$
which proves Proposition~\ref{prop1} in the case where $m=1$, and our induction hypothesis in the basis case.

Assuming in the following that $m \geq 2$, we shall now work on the term 
$$w_1\Big(\textrm{Re }q\big((\textrm{Im }F)^{m-1}X\big)\langle X \rangle^{-\frac{2(2m-1)}{2m+1}}\Big)\langle X \rangle^{\frac{2}{2m+1}},$$
appearing in (\ref{giu35}).
By considering some constants $\Lambda_j \geq 1$, for $0 \leq j \leq m-2$, whose values will be successively chosen in the following, we shall prove that one can write that for all $X \in \rr^{2n}$,
\inc\begin{multline*}\label{giu38}
w_1\left(\frac{\textrm{Re }q\big((\textrm{Im }F)^{m-1}X\big)}{\langle X \rangle^{\frac{2(2m-1)}{2m+1}}}\right) \leq \tilde{W}_0(X)\Psi_0(X) \\ + \sum_{j=1}^{m-2}{\tilde{W}_0(X)
\Big(\prod_{l=1}^{j}W_l(X)\Big)\Psi_j(X)}
+\tilde{W}_{0}(X)\Big(\prod_{l=1}^{m-1}W_l(X)\Big), \num
\end{multline*}
with 
\begin{equation}\label{giu39}\inc
\Psi_j(X)=\psi\left(\frac{\Lambda_j \textrm{Re }q\big((\textrm{Im }F)^{m-j-2}X\big)}{\textrm{Re }q\big((\textrm{Im }F)^{m-j-1}X\big)^{\frac{2m-2j-3}{2m-2j-1}}}\right), \ 0 \leq j \leq m-2, \num
\end{equation} 
\begin{equation}\label{giu40}\inc
W_j(X)=w_{2}\left(\frac{\Lambda_{j-1} \textrm{Re }q\big((\textrm{Im }F)^{m-j-1}X\big)}{\textrm{Re }q\big((\textrm{Im }F)^{m-j}X\big)^{\frac{2m-2j-1}{2m-2j+1}}}\right), \ 1 \leq j \leq m-1,\num
\end{equation}
\begin{equation}\label{giu41}\inc
\tilde{W}_0(X)=w_{1}\left(\frac{\textrm{Re }q\big((\textrm{Im }F)^{m-1}X\big)}{\langle X \rangle^{\frac{2(2m-1)}{2m+1}}}\right),\num
\end{equation}
where $\psi$ is the $C^{\infty}(\rr,[0,1])$ function defined in (\ref{giu13}), and $w_2$ is a $C^{\infty}(\rr,[0,1])$ function satisfying similar properties as the function defined in (\ref{giu15}), with possibly different positive numerical values for its support localization, in order to have that 
\begin{equation}\label{bell1}\inc
\textrm{supp }\psi' \subset \big\{w_2=1\big\} \textrm{ and } \textrm{supp }w_2' \subset \big\{\psi=1\big\}. \num
\end{equation} 
In order to check (\ref{giu38}), we begin by noticing from (\ref{giu15}), (\ref{giu40}) and (\ref{giu41}) that for $0 \leq j \leq m-1$, 
\inc\begin{multline*}\label{bell2}
\textrm{Re }q\big((\textrm{Im }F)^{m-j-1}X\big)^{\frac{1}{2m-2j-1}} \gtrsim \textrm{Re }q\big((\textrm{Im }F)^{m-j}X\big)^{\frac{1}{2m-2j+1}} \\ \gtrsim ... \gtrsim 
\textrm{Re }q\big((\textrm{Im }F)^{m-1}X\big)^{\frac{1}{2m-1}} \gtrsim \langle X \rangle^{\frac{2}{2m+1}}, \num 
\end{multline*}
on the support of the function
$$\textrm{supp} \Big(\tilde{W}_0\prod_{l=1}^{j}W_l\Big), \textrm{ if } 1 \leq j \leq m-1, \textrm{ or, } \textrm{supp }\tilde{W}_0,\textrm{ if } j=0.$$ 
Notice that the constants in the estimates (\ref{bell2}) only depend on the values of the parameters $\Lambda_0$,...,$\Lambda_{j-1}$ but not on $\Lambda_l$, when $l \geq j$.
This shows that the functions 
$$\Psi_0; \ \Big(\prod_{l=1}^{j}W_l\Big)\Psi_j, \textrm{ for } 1 \leq j \leq m-2; \textrm{ and }  \prod_{l=1}^{m-1}W_l,$$
are well-defined on the support of the function $\tilde{W}_0$. Now, by noticing from (\ref{giu13}), (\ref{giu15}), (\ref{giu39}), (\ref{giu40}) and (\ref{bell1}) that
\begin{equation}\label{eg1}\inc
1 \leq \Psi_j+W_{j+1},\num
\end{equation} 
on the support of the function
$$\textrm{supp} \Big(\tilde{W}_0\prod_{l=1}^{j}W_l\Big), \textrm{ if } 1 \leq j \leq m-2, \textrm{ or, } \textrm{supp }\tilde{W}_0,\textrm{ if } j=0,$$ 
we deduce the estimate (\ref{giu38}) from a finite iteration by using the following estimates
$$\tilde{W}_0 \leq \tilde{W}_0\Psi_0+\tilde{W}_0W_1$$
and
$$\tilde{W}_{0}\Big(\prod_{l=1}^{j}W_l\Big) \leq \tilde{W}_{0}\Big(\prod_{l=1}^{j}W_l\Big)\Psi_j+\tilde{W}_{0}\Big(\prod_{l=1}^{j+1}W_l\Big),$$
for any $1 \leq j \leq m-2$. One can also notice that (\ref{eg1}) implies that 
\begin{equation}\label{eg2}\inc
1 \leq \Psi_j+\sum_{k=j+1}^{m-2}{\Big(\prod_{l=j+1}^kW_l\Big)\Psi_k}+\prod_{l=j+1}^{m-1}W_l,\num
\end{equation}
on the support of the function
$$\textrm{supp} \Big(\tilde{W}_0\prod_{l=1}^{j}W_l\Big), \textrm{ if } 1 \leq j \leq m-2, \textrm{ or, } \textrm{supp }\tilde{W}_0,\textrm{ if } j=0.$$ 
Since $\textrm{Re }q \geq 0$, we then get from (\ref{bell2}) that 
\begin{equation}\label{giu42}\inc 
\forall X \in \rr^{2n}, \ \tilde{W}_{0}(X)\Big(\prod_{l=1}^{m-1}W_l(X)\Big)\langle X \rangle^{\frac{2}{2m+1}} \leq \tilde{a}_{\Lambda_0,...,\Lambda_{m-2}}  \textrm{Re }q(X), \num
\end{equation}
where $\tilde{a}_{\Lambda_0,...,\Lambda_{m-2}}$ is a positive constant whose value depend on the parameters 
$$(\Lambda_l)_{0 \leq l \leq m-2}.$$
We define
\begin{equation}\label{giu44}\inc
p_j(X)=\tilde{W}_0(X)\Big(\prod_{l=1}^{j}{W_l(X)}\Big)\Psi_j(X)\frac{r_{m-j-1}(X)}{\textrm{Re }q\big((\textrm{Im }F )^{m-j-1}X\big)^{\frac{2m-2j-2}{2m-2j-1}}},  \num
\end{equation}
for $1 \leq j \leq m-2$, and
\begin{equation}\label{giu45}\inc
p_0(X)=\tilde{W}_0(X)\Psi_0(X)\frac{r_{m-1}(X)}{\textrm{Re }q\big((\textrm{Im }F )^{m-1}X\big)^{\frac{2m-2}{2m-1}}}, \num
\end{equation}
where the quadratic forms $r_{k}$ are defined in (\ref{giu11}). We get from (\ref{giu13}), (\ref{giu15}), (\ref{giu39}), (\ref{giu40}), (\ref{giu41}), (\ref{bell2}), Lemma~\ref{lem3},  Lemma~\ref{lem2.1}, Lemma~\ref{lem2.2} and Lemma~\ref{lem2.3} that 
\begin{equation}\label{giu47}\inc
p_j \in S\big(1,\langle X \rangle^{-\frac{2(2m-2j-3)}{2m+1}}dX^2\big). \num
\end{equation}
for any $0 \leq j \leq m-2$.

We shall now study the Poisson brackets $H_{\textrm{Im}q}\ p_j$. In doing so, we begin by writing that
\inc\begin{align*}\label{giu48}
& \ H_{\textrm{Im}q} \ p_j(X)=\big(H_{\textrm{Im}q}\tilde{W}_0\big)(X)\Big(\prod_{l=1}^{j}{W_l(X)}\Big)\Psi_j(X) \frac{r_{m-j-1}(X)}{\textrm{Re }q\big((\textrm{Im }F )^{m-j-1}X\big)^{\frac{2m-2j-2}{2m-2j-1}}} \num \\
+ & \ \tilde{W}_0(X)\Big(\prod_{l=1}^{j}{W_l(X)}\Big)\big(H_{\textrm{Im}q}\Psi_j\big)(X) \frac{r_{m-j-1}(X)}{\textrm{Re }q\big((\textrm{Im }F )^{m-j-1}X\big)^{\frac{2m-2j-2}{2m-2j-1}}} \\
+ & \  \tilde{W}_0(X)\Big(\prod_{l=1}^{j}{W_l(X)}\Big)\Psi_j(X)H_{\textrm{Im}q}\left(\textrm{Re }q\big((\textrm{Im }F )^{m-j-1}X\big)^{-\frac{2m-2j-2}{2m-2j-1}}\right)r_{m-j-1}(X) \\
+ & \ \tilde{W}_0(X)\Big(\prod_{l=1}^{j}{W_l(X)}\Big)\Psi_j(X)\frac{H_{\textrm{Im}q}\ r_{m-j-1}(X)}{\textrm{Re }q\big((\textrm{Im }F )^{m-j-1}X\big)^{\frac{2m-2j-2}{2m-2j-1}}}\\
+ & \ \sum_{l=1}^j\tilde{W}_0(X)\big(H_{\textrm{Im}q}W_l\big)(X)\Big(\prod_{k=1 \atop k \neq l}^{j}{W_k(X)}\Big)\Psi_j(X) \frac{r_{m-j-1}(X)}{\textrm{Re }q\big((\textrm{Im }F )^{m-j-1}X\big)^{\frac{2m-2j-2}{2m-2j-1}}},
\end{align*}
for $1 \leq j \leq m-2$. We denote by respectively $B_{1,j}$, $B_{2,j}$, $B_{3,j}$, $B_{4,j}$ and $B_{5,j}$ the five terms appearing in the right hand side of (\ref{giu48}). 
We also write in the case where $j=0$,
\inc\begin{align*}\label{giu48bis}
& \ H_{\textrm{Im}q} \ p_0(X)=\big(H_{\textrm{Im}q}\tilde{W}_0\big)(X)\Psi_0(X) \frac{r_{m-1}(X)}{\textrm{Re }q\big((\textrm{Im }F )^{m-1}X\big)^{\frac{2m-2}{2m-1}}} \num \\
+ & \ \tilde{W}_0(X)\big(H_{\textrm{Im}q}\Psi_0\big)(X) \frac{r_{m-1}(X)}{\textrm{Re }q\big((\textrm{Im }F )^{m-1}X\big)^{\frac{2m-2}{2m-1}}} \\
+ & \  \tilde{W}_0(X)\Psi_0(X)H_{\textrm{Im}q}\left(\textrm{Re }q\big((\textrm{Im }F )^{m-1}X\big)^{-\frac{2m-2}{2m-1}}\right)r_{m-1}(X) \\
+ & \ \tilde{W}_0(X)\Psi_0(X)\frac{H_{\textrm{Im}q}\ r_{m-1}(X)}{\textrm{Re }q\big((\textrm{Im }F )^{m-1}X\big)^{\frac{2m-2}{2m-1}}},
\end{align*}
and denote as before by respectively $B_{1,0}$, $B_{2,0}$, $B_{3,0}$ and $B_{4,0}$ the four terms appearing in the right hand side of (\ref{giu48bis}).

Since the constants in the estimates (\ref{bell2}) only depend on the values of the parameters $\Lambda_0$,..., $\Lambda_{j-1}$; but not on $\Lambda_l$, when $l \geq j$; we notice from (\ref{giu38}), (\ref{bell2}) and (\ref{giu42}) that there exist $a_{0}>0$ and some positive constants $a_{j,\Lambda_0,...,\Lambda_{j-1}}$, for $1 \leq j \leq m-1$, whose values with respect to the parameters $(\Lambda_l)_{0 \leq l \leq m-2}$ only depend on $\Lambda_0$,..., $\Lambda_{j-1}$; but not on $\Lambda_l$, when $l \geq j$; such that for any constants $(\alpha_j)_{1 \leq j \leq m-2}$, with $\alpha_j \geq 1$; and $X \in \rr^{2n}$,
\inc\begin{align*}\label{giu100}
 \num & \ w_1\left(\frac{\textrm{Re }q\big((\textrm{Im }F)^{m-1}X\big)}{\langle X \rangle^{\frac{2(2m-1)}{2m+1}}}\right) \langle X \rangle^{\frac{2}{2m+1}} 
\leq a_0 \tilde{W}_0(X)\Psi_0(X) \textrm{Re }q\big((\textrm{Im }F)^{m-1}X\big)^{\frac{1}{2m-1}}\\ 
 & \  + \sum_{j=1}^{m-2}{\alpha_j a_{j,\Lambda_0,...,\Lambda_{j-1}} \tilde{W}_0(X)
\Big(\prod_{l=1}^{j}W_l(X)\Big)\Psi_j(X)}\textrm{Re }q\big((\textrm{Im }F)^{m-j-1}X\big)^{\frac{1}{2m-2j-1}}\\
  & \ +a_{m-1,\Lambda_0,...,\Lambda_{m-2}}\textrm{Re }q(X).
\end{align*}
The positive constant $a_0$ is independent of any of the parameters $(\Lambda_l)_{0 \leq l \leq m-2}$.
Setting 
\begin{equation}\label{giu101}\inc
p=a_0 p_0+\sum_{j=1}^{m-2}{\alpha_j a_{j,\Lambda_0,...,\Lambda_{j-1}} p_j}, \num
\end{equation}
we know from (\ref{giu47}) that 
\begin{equation}\label{giu102}\inc
p \in S\big(1,\langle X \rangle^{-\frac{2}{2m+1}}dX^2\big).\num
\end{equation}
For any $\eps >0$, we shall prove that after a proper choice for the constants $(\Lambda_j)_{0 \leq j \leq m-2}$ and $(\alpha_j)_{1 \leq j \leq m-2}$, with $\Lambda_j \geq1$, $\alpha_j \geq 1$, whose values will depend on $\eps$; one can find a positive constant $c_{12,\eps}>0$ such that for all $X \in \rr^{2n}$,
\inc\begin{equation}\label{giu103}
 c_{12,\eps}\textrm{Re }q(X)+H_{\textrm{Im}q} \ p(X) +\eps \langle X \rangle^{\frac{2}{2m+1}} \geq  w_1\left(\frac{\textrm{Re }q\big((\textrm{Im }F)^{m-1}X\big)}{\langle X \rangle^{\frac{2(2m-1)}{2m+1}}}\right) \langle X \rangle^{\frac{2}{2m+1}}.\num
\end{equation}
Once this estimate proved, Proposition~\ref{prop1} will directly follow from (\ref{giu33}), (\ref{giu35}), (\ref{giu102}) and (\ref{giu103}), if we choose the positive parameter $\eps$ sufficiently small and consider the weight function
$$g=c_{13,\eps}G_m+c_{14,\eps}p,$$ 
after a suitable choice for the positive constants $c_{13,\eps}$ and $c_{14,\eps}$.

Let $\eps>0$, it therefore remains to choose properly these constants $(\Lambda_j)_{0 \leq j \leq m-2}$ and $(\alpha_j)_{1 \leq j \leq m-2}$, with $\Lambda_j \geq1$, $\alpha_j \geq 1$, in order to satisfy (\ref{giu103}).

Recalling from (\ref{kee1}) that 
\inc\begin{multline*}\label{giu104}
\forall \ 0 \leq j \leq m-2, \ H_{\textrm{Im}q}\ r_{m-j-1}(X)=2\textrm{Re }q\big((\textrm{Im }F )^{m-j-1}X\big)\\
+2\textrm{Re }q\big((\textrm{Im }F )^{m-j}X;(\textrm{Im }F )^{m-j-2}X \big),\num
\end{multline*}
one can notice by expanding the term $2a_{m-1,\Lambda_0,...,\Lambda_{m-2}}\textrm{Re }q+H_{\textrm{Im}q}\ p$ by using (\ref{giu48}), (\ref{giu48bis}) and (\ref{giu101}) that the terms in
$$a_0 B_{4,0}+\sum_{j=1}^{m-2}{\alpha_j a_{j,\Lambda_0,...,\Lambda_{j-1}} B_{4,j}},$$
produced by the terms associated to $2\textrm{Re }q\big((\textrm{Im }F )^{m-j-1}X\big)$ while using (\ref{giu104}), give exactly two times the term 
\inc\begin{align*}\label{giu105}
& \ a_0\tilde{W}_0(X)\Psi_0(X) \textrm{Re }q\big((\textrm{Im }F)^{m-1}X\big)^{\frac{1}{2m-1}} \num \\ 
+ & \  \sum_{j=1}^{m-2}{\alpha_j a_{j,\Lambda_0,...,\Lambda_{j-1}} \tilde{W}_0(X)
\Big(\prod_{l=1}^{j}W_l(X)\Big)\Psi_j(X)}\textrm{Re }q\big((\textrm{Im }F)^{m-j-1}X\big)^{\frac{1}{2m-2j-1}}\\
+ & \ a_{m-1,\Lambda_0,...,\Lambda_{m-2}}\textrm{Re }q(X),
\end{align*}
for which we have the estimate (\ref{giu100}). To prove the estimate (\ref{giu103}), it will therefore be sufficient to check that all the other terms appearing in (\ref{giu48}) and (\ref{giu48bis}) can also be all absorbed in the term (\ref{giu105}) after a proper choice for the constants $(\Lambda_j)_{0 \leq j \leq m-2}$ and $(\alpha_j)_{1 \leq j \leq m-2}$; at the exception of a remainder term in 
$$\eps \langle X \rangle^{\frac{2}{2m+1}}.$$ 
We shall choose these constants in the following order $\Lambda_0$, $\alpha_1$, $\Lambda_1$, $\alpha_2$, ...., $\alpha_{m-2}$ and $\Lambda_{m-2}$. 

We successively study the remaining terms in (\ref{giu48}) and (\ref{giu48bis}), by increasing value of the integer $0 \leq j \leq m-2$. We first notice from (\ref{giu13}), (\ref{giu15}), (\ref{giu39}), (\ref{giu41}), (\ref{giu48bis}), Lemma~\ref{lem2.331} and Lemma~\ref{lem2.435} that one can choose the first constant $\Lambda_0 \geq 1$ such that for all $X \in \rr^{2n}$,
\begin{equation}\label{giu106}\inc
a_0 |B_{1,0}(X)| \lesssim \Lambda_0^{-\frac{1}{2}} \langle X \rangle^{\frac{2}{2m+1}} \leq \frac{\eps}{m-1} \langle X \rangle^{\frac{2}{2m+1}}.\num
\end{equation}
By noticing from (\ref{giu15}) and (\ref{giu41}) that the estimates
\begin{equation}\label{giu107}\inc
\textrm{Re }q\big((\textrm{Im }F )^{m}X\big) \lesssim \langle X \rangle^2 \lesssim \textrm{Re }q\big((\textrm{Im }F )^{m-1}X\big)^{\frac{2m+1}{2m-1}},\num
\end{equation}
are fulfilled on the support of the function $\tilde{W}_0$, we deduce from (\ref{giu13}), (\ref{giu15}), (\ref{giu39}), (\ref{giu41}), (\ref{giu48bis}), Lemma~\ref{lem2.331}, Lemma~\ref{lem2.4} and Lemma~\ref{lem2.41} that the modulus of the term $B_{3,0}$ and the second term in $B_{4,0}$ associated to 
$$2\textrm{Re }q\big((\textrm{Im }F )^{m}X;(\textrm{Im }F )^{m-2}X \big),$$ 
while using (\ref{giu104}), that we denote here $\tilde{B}_{4,0}$, can both be estimated as
$$a_0|B_{3,0}(X)|+a_0|\tilde{B}_{4,0}(X)| \lesssim \Lambda_0^{-\frac{1}{2}}\tilde{W}_0(X)\Psi_0(X)\textrm{Re }q\big((\textrm{Im }F )^{m-1}X\big)^{\frac{1}{2m-1}},$$
for all $X \in \rr^{2n}$. By possibly increasing sufficiently the value of the constant $\Lambda_0$ which is of course possible while keeping (\ref{giu106}), one can control these terms with the \og good \fg \ term (\ref{giu105}). The value of the constant $\Lambda_0$ is now definitively fixed. In (\ref{giu48bis}), it only remains to study the term $B_{2,0}$.

About this term, we deduce from (\ref{giu13}), (\ref{giu15}), (\ref{giu39}), (\ref{giu41}), (\ref{giu48bis}), Lemma \ref{lem2.331} and Lemma~\ref{lem2.42} that for all $X \in \rr^{2n}$,
\begin{equation}\label{eg3}\inc
a_0|B_{2,0}(X)| \lesssim \tilde{W}_0(X)W_1(X)\textrm{Re }q\big((\textrm{Im }F)^{m-1}X\big)^{\frac{1}{2m-1}}.\num
\end{equation}
By using now (\ref{bell2}) and (\ref{eg2}) with $j=1$, we obtain that for all $X \in \rr^{2n}$,
\begin{multline*}
a_0|B_{2,0}(X)| \leq c_{m-1,\Lambda_0,...,\Lambda_{m-2}}\tilde{W}_0(X)\Big(\prod_{l=1}^{m-1}W_l(X)\Big)\textrm{Re }q(X) \\
+ \sum_{j=1}^{m-2}{c_{j,\Lambda_0,...,\Lambda_{j-1}}\tilde{W}_0(X)\Big(\prod_{l=1}^jW_l(X)\Big)\Psi_j(X)\textrm{Re }q\big((\textrm{Im }F)^{m-j-1}X\big)^{\frac{1}{2m-2j-1}}},
\end{multline*}
which implies that 
\inc\begin{multline*}\label{eg4}
a_0|B_{2,0}(X)| \leq c_{m-1,\Lambda_0,...,\Lambda_{m-2}}\textrm{Re }q(X)\\
+ \sum_{j=1}^{m-2}{c_{j,\Lambda_0,...,\Lambda_{j-1}}\tilde{W}_0(X)\Big(\prod_{l=1}^jW_l(X)\Big)\Psi_j(X)\textrm{Re }q\big((\textrm{Im }F)^{m-j-1}X\big)^{\frac{1}{2m-2j-1}}},\num
\end{multline*}
where the quantities $c_{j,\Lambda_0,...,\Lambda_{j-1}}$ stand for positive constants whose values depend on $\Lambda_0$,..., $\Lambda_{j-1}$, but not on 
$(\Lambda_k)_{j \leq k \leq m-2}$ and $(\alpha_k)_{1 \leq k \leq m-2},$ 
according to the remark done after (\ref{bell2}). One can therefore choose the constant $\alpha_1 \geq 1$ in (\ref{giu101}) sufficiently large in order to absorb the term of the index $j=1$ in the sum appearing in the right hand side of the estimate (\ref{eg4}) by the term of same index in the \og good\fg \ term (\ref{giu105}). This is possible since the constants $a_{1,\Lambda_0}$ and 
$c_{1,\Lambda_0}$ are now fixed after our choice of the parameter $\Lambda_0$.

This ends our step index $j=0$ in which we have chosen the values for the two constants $\Lambda_0$ and $\alpha_1 \geq 1$. We shall now explain how to choose the remaining constants
$(\Lambda_j)_{1 \leq j \leq m-2}$ and $(\alpha_j)_{2 \leq j \leq m-2}$ in (\ref{giu101}) in order to satisfy (\ref{giu103}). This choice will also determine the values of the constants $(a_{j,\Lambda_0,...,\Lambda_{j-1}})_{1 \leq j \leq m-2}$ appearing in (\ref{giu101}). After this step index $j=0$, we have managed to absorb all the terms appearing in (\ref{giu48bis}) in the \og good \fg \ term (\ref{giu105}) at the exception of a remainder coming from (\ref{giu106}) and (\ref{eg4}),
\begin{multline*}
\sum_{j=2}^{m-2}{c_{j,\Lambda_0,...,\Lambda_{j-1}}\tilde{W}_0(X)\Big(\prod_{l=1}^jW_l(X)\Big)\Psi_j(X)\textrm{Re }q\big((\textrm{Im }F)^{m-j-1}X\big)^{\frac{1}{2m-2j-1}}}\\ 
+\frac{\eps}{m-1}\langle X\rangle^{\frac{2}{2m+1}},
\end{multline*}
where one recall that the positive constants $c_{j,\Lambda_0,...,\Lambda_{j-1}}$ only depend on $\Lambda_0$,...,$\Lambda_{j-1}$, but not on 
$(\Lambda_k)_{j \leq k \leq m-2}$ and $(\alpha_k)_{1 \leq k \leq m-2}.$

We proceed in the following by finite induction and assume that, at the beginning of the step index $k$, with $1 \leq k \leq m-2$, we have already chosen the values for the constants 
$(\Lambda_j)_{0 \leq j \leq k-1}$ and $(\alpha_j)_{1 \leq j \leq k}$ in (\ref{giu101}); and that these choices have allowed to absorb all the terms appearing in the right hand side of (\ref{giu48bis}) and (\ref{giu48}), when $1 \leq j \leq k-1$, in the \og good \fg \ term (\ref{giu105}) at the exception of a remainder term
\inc\begin{multline*}\label{eg5}
\frac{k}{m-1}\eps \langle X \rangle^{\frac{2}{2m+1}}+\\
\sum_{j=k+1}^{m-2}\tilde{c}_{j,\Lambda_0,...,\Lambda_{j-1}, \alpha_1,...,\alpha_{k-1}}\tilde{W}_0(X)\Big(\prod_{l=1}^jW_l(X)\Big)\Psi_j(X)\textrm{Re }q\big((\textrm{Im }F)^{m-j-1}X\big)^{\frac{1}{2m-2j-1}},\num
\end{multline*} 
where the quantities $\tilde{c}_{j,\Lambda_0,...,\Lambda_{j-1}, \alpha_1,...,\alpha_{k-1}}$ stand for positive constants whose values only depend on $\Lambda_0$,..., $\Lambda_{j-1}$, $\alpha_1$,..., $\alpha_{k-1}$; but not on $(\Lambda_l)_{j \leq l \leq m-2}$ and $(\alpha_l)_{k \leq l \leq m-2}$.

We shall now explain how to choose the constants $\Lambda_k$ and; $\alpha_{k+1}$, when $k \leq m-3$; in this step index $k$ in order to absorb the terms appearing in the right hand side of (\ref{giu48}), when $j=k$, at the exception of a remainder term of the type (\ref{eg5}) where $k$ will be replaced by $k+1$;  in the \og good \fg \ term (\ref{giu105}). Since the constants $(\Lambda_j)_{0 \leq j \leq k-1}$ and $(\alpha_j)_{1 \leq j \leq k}$ have already been chosen, we shall only underline in the following the dependence of our estimates with respect to the other parameters 
$(\Lambda_j)_{k \leq j \leq m-2}$ and $(\alpha_j)_{k+1 \leq j \leq m-2}$, whose values remain to be chosen.

We notice from (\ref{giu13}), (\ref{giu15}), (\ref{giu39}), (\ref{giu40}), (\ref{giu41}), (\ref{bell2}), (\ref{giu48}), Lemma~\ref{lem2.331} and Lemma~\ref{lem2.435} that one can assume by choosing the constant $\Lambda_k \geq 1$ sufficiently large that for all $X \in \rr^{2n}$,
\begin{equation}\label{eg6}\inc
\alpha_k a_{k,\Lambda_0,...,\Lambda_{k-1}} |B_{1,k}(X)| \lesssim \Lambda_k^{-\frac{1}{2}} \langle X \rangle^{\frac{2}{2m+1}} \leq \frac{\eps}{m-1} \langle X \rangle^{\frac{2}{2m+1}},\num
\end{equation}
since the constants $\alpha_k$, $\Lambda_0$,....,$\Lambda_{k-1}$ have already been fixed.
Next, we deduce from (\ref{giu13}), (\ref{giu15}), (\ref{giu39}), (\ref{giu40}), (\ref{giu41}), (\ref{bell2}), (\ref{giu48}), Lemma~\ref{lem2.331}, Lemma~\ref{lem2.4} and Lemma~\ref{lem2.41} that the modulus of the term $B_{3,k}$ and the second term in $B_{4,k}$ associated to 
$$2\textrm{Re }q\big((\textrm{Im }F )^{m-k}X;(\textrm{Im }F )^{m-k-2}X \big),$$ 
while using (\ref{giu104}), that we denote here $\tilde{B}_{4,k}$, can both be estimated as
\begin{multline*}
\alpha_k a_{k,\Lambda_0,...,\Lambda_{k-1}}|B_{3,k}(X)|+\alpha_k a_{k,\Lambda_0,...,\Lambda_{k-1}}|\tilde{B}_{4,k}(X)| \\ 
\lesssim \Lambda_k^{-\frac{1}{2}}\tilde{W}_0(X)\Big(\prod_{l=1}^kW_l(X)\Big)\Psi_k(X)\textrm{Re }q\big((\textrm{Im }F )^{m-k-1}X\big)^{\frac{1}{2m-2k-1}},
\end{multline*}
for all $X \in \rr^{2n}$. By possibly increasing sufficiently the value of the constant $\Lambda_k$ which is of course possible while keeping (\ref{eg6}), one can control these terms with the \og good \fg \ term (\ref{giu105}).

For $1 \leq l \leq k$, we shall now study the term
$$B_{5,k,l}(X)=\tilde{W}_0(X)\big(H_{\textrm{Im}q}W_l\big)(X)\Big(\prod_{j=1 \atop j \neq l}^{k}{W_j(X)}\Big)\Psi_k(X) \frac{r_{m-k-1}(X)}{\textrm{Re }q\big((\textrm{Im }F )^{m-k-1}X\big)^{\frac{2m-2k-2}{2m-2k-1}}},$$
appearing in the term $B_{5,k}$ in (\ref{giu48}). 
By noticing that 
$$\textrm{Re }q\big((\textrm{Im }F)^{m-l-2}X\big) \sim \Lambda_l^{-1}\textrm{Re }q\big((\textrm{Im }F)^{m-l-1}X\big)^{\frac{2m-2l-3}{2m-2l-1}},$$
on the support of the function $H_{\textrm{Im}q}W_{l+1}$, it follows from (\ref{giu13}), (\ref{giu15}), (\ref{giu39}), (\ref{giu40}), (\ref{giu41}), (\ref{bell2}), (\ref{giu107}), Lemma \ref{lem2.331}
and Lemma~\ref{lem2.43} that for all $X \in \rr^{2n}$,
$$\alpha_k a_{k,\Lambda_0,...,\Lambda_{k-1}}|B_{5,k,1}(X)| \lesssim \Lambda_k^{-\frac{1}{2}}\tilde{W}_0(X)\Psi_{0}(X)\textrm{Re }q\big((\textrm{Im }F )^{m-1}X\big)^{\frac{1}{2m-1}}$$
and
\begin{multline*}
\alpha_k a_{k,\Lambda_0,...,\Lambda_{k-1}}|B_{5,k,l}(X)| \\ 
\lesssim \Lambda_k^{-\frac{1}{2}}\tilde{W}_0(X)\Big(\prod_{j=1}^{l-1}{W_j(X)}\Big)\Psi_{l-1}(X)\textrm{Re }q\big((\textrm{Im }F )^{m-l}X\big)^{\frac{1}{2m-2l+1}},
\end{multline*}
when $l \geq 2$. By possibly increasing again the value of the constant $\Lambda_k$, one can therefore control the term $\alpha _k a_{k,\Lambda_0,...,\Lambda_{k-1}}B_{5,k}$ with the \og good \fg \ term (\ref{giu105}).
The value of the constant $\Lambda_k$ is now definitively fixed.

About the term $B_{2,k}$, we deduce from (\ref{giu13}), (\ref{giu15}), (\ref{giu39}), (\ref{giu40}), (\ref{giu41}), (\ref{bell2}), (\ref{giu48}), Lemma \ref{lem2.331} and Lemma~\ref{lem2.42} that for all $X \in \rr^{2n}$,
\inc\begin{multline*}\label{eg7}
\alpha_k a_{k,\Lambda_0,...,\Lambda_{k-1}}|B_{2,k}(X)| \\ \lesssim \tilde{W}_0(X)\Big(\prod_{l=1}^{k+1}W_l(X)\Big)\textrm{Re }q\big((\textrm{Im }F)^{m-k-1}X\big)^{\frac{1}{2m-2k-1}}.\num
\end{multline*}
By distinguishing two cases, we first assume in the following that $k \leq m-3$. In this case,
by using (\ref{bell2}) and (\ref{eg2}) with $j=k+1$, we obtain that for all $X \in \rr^{2n}$,
\begin{multline*}
\alpha_k a_{k,\Lambda_0,...,\Lambda_{k-1}}|B_{2,k}(X)| \leq  c_{m-1,\Lambda_0,...,\Lambda_{m-2},\alpha_1,...,\alpha_k}'\tilde{W}_0(X)\Big(\prod_{l=1}^{m-1}W_l(X)\Big)\textrm{Re }q(X)\\
+\sum_{j=k+1}^{m-2}{c_{j,\Lambda_0,...,\Lambda_{j-1},\alpha_1,...,\alpha_k}'\tilde{W}_0(X)\Big(\prod_{l=1}^jW_l(X)\Big)\Psi_j(X)\textrm{Re }q\big((\textrm{Im }F)^{m-j-1}X\big)^{\frac{1}{2m-2j-1}}},
\end{multline*}
which implies that 
\inc\begin{multline*}\label{eg8}
\alpha_k a_{k,\Lambda_0,...,\Lambda_{k-1}}|B_{2,k}(X)| \leq  c_{m-1,\Lambda_0,...,\Lambda_{m-2},\alpha_1,...,\alpha_k}'\textrm{Re }q(X)\\
+\sum_{j=k+1}^{m-2}{c_{j,\Lambda_0,...,\Lambda_{j-1},\alpha_1,...,\alpha_k}'\tilde{W}_0(X)\Big(\prod_{l=1}^jW_l(X)\Big)\Psi_j(X)\textrm{Re }q\big((\textrm{Im }F)^{m-j-1}X\big)^{\frac{1}{2m-2j-1}}},\num
\end{multline*}
where the quantities $c_{j,\Lambda_0,...,\Lambda_{j-1},\alpha_1,...,\alpha_k}'$ stand for positive constants whose values only depend on $\Lambda_0$,..., $\Lambda_{j-1}$, $\alpha_1$,..., $\alpha_{k}$, but not on $(\Lambda_l)_{j \leq l \leq m-2}$ and $(\alpha_l)_{k+1 \leq l \leq m-2}$. Indeed, we recall that the constants appearing in the estimates (\ref{bell2}) only depend on the values of the parameters $\Lambda_0$,..., $\Lambda_{j-1}$; but not on $(\Lambda_l)_{j \leq l \leq m-2}$ and $(\alpha_l)_{1 \leq l \leq m-2}$. One can therefore choose the constant $\alpha_{k+1} \geq 1$ in (\ref{giu101}) sufficiently large in order to absorb the term of index $j=k+1$ in the sum (\ref{eg5});
and the term of index $j=k+1$ in the sum appearing in the right hand side of the estimate (\ref{eg8}), by the term of same index in the \og good\fg \ term (\ref{giu105}).

When $k=m-2$ and taking $\Lambda_{m-2}=1$, it follows from (\ref{bell2}), used with $j=m-1$, and (\ref{eg7}) that for all $X \in \rr^{2n}$,
\inc\begin{align*}\label{eg9}
\alpha_{m-2} a_{m-2,\Lambda_0,...,\Lambda_{m-3}} |B_{2,m-2}(X)| \lesssim & \ \tilde{W}_0(X)\Big(\prod_{l=1}^{m-1}W_l(X)\Big)\textrm{Re }q(\textrm{Im }FX)^{\frac{1}{3}} \num \\
\lesssim & \  \textrm{Re }q(X).
\end{align*}
This process allows us to achieve the construction of the weight function $p$ satisfying (\ref{giu103}), which ends the proof of (\ref{giu103}). This also ends the proof of Proposition~\ref{prop1}.~$\Box$

\section{Appendix}\label{appendix}
\init

\subsection{Wick calculus}\label{wick} The purpose of this section is to recall the definition and basic properties of the Wick quantization that we need for the proof of Theorem~\ref{theorem1}. We follow here the presentation of the Wick quantization given by N.~Lerner in \cite{lerner} and refer the reader to his work for the proofs of the results recalled below.

The main property of the Wick quantization is its property of positivity, i.e., that non-negative Hamiltonians define non-negative operators
$$a \geq 0 \Rightarrow a^{\textrm{Wick}} \geq 0.$$
We recall that this is not the case for the Weyl quantization and refer to \cite{lerner} for an example of non-negative Hamiltonian defining an operator which is not non-negative.

Before defining properly the Wick quantization, we first need to recall the definition of the wave packets transform of a function $u \in \mathcal{S}(\rr^n)$, 
$$Wu(y,\eta)=(u,\varphi_{y,\eta})_{L^2(\rr^n)}=2^{n/4}\int_{\rr^n}{u(x)e^{- \pi (x-y)^2}e^{-2i \pi(x-y).\eta}dx}, \ (y,\eta) \in \rr^{2n}.$$
where 
$$\varphi_{y,\eta}(x)=2^{n/4}e^{- \pi (x-y)^2}e^{2i \pi (x-y).\eta}, \ x \in \mathbb{R}^n,$$
and $x^2=x_1^2+...+x_n^2$. With this definition, one can check (see Lemma 2.1 in \cite{lerner}) that 
the mapping $u \mapsto Wu$ is continuous from $\mathcal{S}(\rr^n)$ to $\mathcal{S}(\rr^{2n})$, isometric from $L^{2}(\rr^n)$ to $L^2(\rr^{2n})$ and that we have the
reconstruction formula
\begin{equation}\label{lay0.1}\inc
\forall u \in \mathcal{S}(\rr^n), \forall x \in \rr^n, \ u(x)=\int_{\rr^{2n}}{Wu(y,\eta)\varphi_{y,\eta}(x)dyd\eta}.\num
\end{equation}
By denoting by $\Sigma_Y$ the operator defined in the Weyl quantization by the symbol 
$$p_Y(X)=2^n e^{-2\pi|X-Y|^2}, \ Y=(y,\eta) \in \rr^{2n},$$
which is a rank-one orthogonal projection,
$$\big{(}\Sigma_Y u\big{)}(x)=Wu(Y)\varphi_Y(x)=(u,\varphi_Y)_{L^2(\rr^n)}\varphi_Y(x),$$
we define the Wick quantization of any $L^{\infty}(\rr^{2n})$  symbol $a$ as
\begin{equation}\label{lay0.2}\inc
a^{\textrm{Wick}}=\int_{\rr^{2n}}{a(Y)\Sigma_Y dY}.\num
\end{equation}
More generally, one can extend this definition when the symbol $a$ belongs to $\mathcal{S}'(\rr^{2n})$ by defining the operator $a^{\textrm{Wick}}$ for any $u$ and $v$ in $\mathcal{S}(\rr^{n})$ by
$$<a^{\textrm{Wick}}u,\overline{v}>_{\mathcal{S}'(\rr^{n}),\mathcal{S}(\rr^{n})}=<a(Y),(\Sigma_Yu,v)_{L^2(\rr^n)}>_{\mathcal{S}'(\rr^{2n}),\mathcal{S}(\rr^{2n})},$$
where $<\textrm{\textperiodcentered},\textrm{\textperiodcentered}>_{\mathcal{S}'(\rr^n),\mathcal{S}(\rr^n)}$ denotes the duality bracket between the
spaces $\mathcal{S}'(\rr^n)$ and $\mathcal{S}(\rr^n)$. The Wick quantization is a positive quantization
\begin{equation}\label{lay0.5}\inc
a \geq 0 \Rightarrow a^{\textrm{Wick}} \geq 0. \num
\end{equation}
In particular, real Hamiltonians get quantized in this quantization by formally self-adjoint operators and one has (see Proposition 3.2 in \cite{lerner}) that $L^{\infty}(\rr^{2n})$ symbols define bounded operators on $L^2(\rr^n)$ such that 
\begin{equation}\label{lay0}\inc
\|a^{\textrm{Wick}}\|_{\mathcal{L}(L^2(\rr^n))} \leq \|a\|_{L^{\infty}(\rr^{2n})}.\num
\end{equation}
According to Proposition~3.3 in~\cite{lerner}, the Wick and Weyl quantizations of a symbol $a$ are linked by the following identities
\begin{equation}\label{lay1bis}\inc
a^{\textrm{Wick}}=\tilde{a}^w,\num
\end{equation}
with
\begin{equation}\label{lay2bis}\inc
\tilde{a}(X)=\int_{\rr^{2n}}{a(X+Y)e^{-2\pi |Y|^2}2^ndY}, \ X \in \rr^{2n},\num
\end{equation}
and
\begin{equation}\label{lay1}\inc
a^{\textrm{Wick}}=a^w+r(a)^w,\num
\end{equation}
where $r(a)$ stands for the symbol
\begin{equation}\label{lay2}\inc
r(a)(X)=\int_0^1\int_{\rr^{2n}}{(1-\theta)a''(X+\theta Y)Y^2e^{-2\pi |Y|^2}2^ndYd\theta}, \ X \in \rr^{2n},\num
\end{equation}
if we use here the normalization chosen in \cite{lerner} for the Weyl quantization
\begin{equation}\label{lay3}\inc
(a^wu)(x)=\int_{\rr^{2n}}{e^{2i\pi(x-y).\xi}a\Big(\frac{x+y}{2},\xi\Big)u(y)dyd\xi},\num
\end{equation}
which differs from the one chosen in this paper. Because of this difference in normalizations, certain constant factors will naturally appear in the core of the proof of Theorem~\ref{theorem1}
while using certain formulas of Section~\ref{wick}, but these are minor adaptations.  
We also recall the following composition formula obtained in the proof of Proposition~3.4 in~\cite{lerner},
\begin{equation}\label{lay4}\inc
a^{\textrm{Wick}} b^{\textrm{Wick}} =\Big{[}ab-\frac{1}{4 \pi} a'.b'+\frac{1}{4i \pi}\{a,b\} \Big{]}^{\textrm{Wick}}+S, \num
\end{equation}
with $\|S\|_{\mathcal{L}(L^2(\rr^n))} \leq d_n \|a\|_{L^{\infty}}\gamma_{2}(b),$
when $a \in L^{\infty}(\rr^{2n})$ and $b$ is a smooth symbol satisfying
$$\gamma_2(b)=\sup_{X \in \rr^{2n}, \atop T \in \rr^{2n}, |T|=1}|b^{(2)}(X)T^2| < +\infty.$$ 
The term $d_n$ appearing in the previous estimate stands for a positive constant depending only on the dimension $n$, and the notation $\{a,b\}$ denotes the Poisson bracket
$$\{a,b\}=\frac{\partial a}{\partial \xi}.\frac{\partial b}{\partial x}-\frac{\partial a}{\partial x}.\frac{\partial b}{\partial \xi}.$$

\subsection{Some technical lemmas}

This second part of the appendix is devoted to the proofs of several technical lemmas.

\bigskip

\begin{lemma}\label{lem2}
If $l_1$, $l_2 \in \nn$ then
\inc\begin{multline*}\label{mari5}
H_{\emph{\textrm{Im}}q} \ \emph{\textrm{Re }}q\big((\emph{\textrm{Im }}F)^{l_1}X;(\emph{\textrm{Im }}F)^{l_2}X\big)=2\emph{\textrm{Re }}q\big((\emph{\textrm{Im }}F)^{l_1+1}X;(\emph{\textrm{Im }}F)^{l_2}X\big)
\\+2\emph{\textrm{Re }}q\big((\emph{\textrm{Im }}F)^{l_1}X;(\emph{\textrm{Im }}F)^{l_2+1}X\big), \num
\end{multline*}
where $\emph{\textrm{Re }}q(X;Y)$ stands for the polarized form associated to the quadratic form $\emph{\textrm{Re }}q$.
\end{lemma}

\bigskip

\noindent
\textit{Proof of Lemma~\ref{lem2}}. 
We begin by noticing from (\ref{10}) and the skew-symmetry property of Hamilton maps (\ref{12}) that the Hamilton map of the quadratic form
$$\tilde{r}(X)=\textrm{Re }q\big((\textrm{Im }F)^{l_1}X;(\textrm{Im }F)^{l_2}X\big),$$ 
is given by 
\begin{equation}\label{lol10}\inc
\tilde{F}=\frac{1}{2}\big((-1)^{l_1}(\textrm{Im }F)^{l_1}\textrm{Re }F(\textrm{Im }F)^{l_2}+(-1)^{l_2}(\textrm{Im }F)^{l_2}\textrm{Re }F(\textrm{Im }F)^{l_1}\big),\num
\end{equation} 
since for any $l_1$, $l_2 \in \nn$,
\inc\begin{align*}\label{mari6}
(-1)^{l_1}\sigma\big(X,(\textrm{Im }F)^{l_1}\textrm{Re }F(\textrm{Im }F)^{l_2}X\big)= & \ \sigma\big((\textrm{Im }F)^{l_1}X,\textrm{Re }F(\textrm{Im }F)^{l_2}X\big) \num \\
= &\ \textrm{Re }q\big((\textrm{Im }F)^{l_1}X,(\textrm{Im }F)^{l_2}X\big).
\end{align*}
Then, a direct computation (see Lemma~2 in~\cite{mz}) shows that the Hamiton map of the quadratic form 
$$H_{\textrm{Im}q} \ \tilde{r}=\big\{\textrm{Im }q,\tilde{r}\big\}=\frac{\partial \textrm{Im }q}{\partial \xi}.\frac{\partial \tilde{r}}{\partial x}-\frac{\partial \textrm{Im } q}{\partial x}.\frac{\partial \tilde{r}}{\partial \xi},$$
is given by the commutator $-2[\textrm{Im }F,\tilde{F}]$, that is,
$$H_{\textrm{Im}q} \ \tilde{r}(X)=-2\sigma\big(X,[\textrm{Im }F,\tilde{F}]X\big).$$
A computation using (\ref{mari6}) then allows to directly get (\ref{mari5}).~$\Box$

\bigskip

\begin{lemma}\label{lem3}
Consider a $C^{\infty}(\rr)$ function $f$ such that 
$$f \in L^{\infty}(\rr) \textrm{ and } \exists c_1,c_2>0, \ \emph{\textrm{supp }}f' \subset \big\{x \in \rr : c_1 \leq |x| \leq c_2\big\},$$
then for all $0<\alpha \leq 1$ and $k \in \nn$,
\begin{equation}\label{giu16}\inc
f\big(\emph{\textrm{Re }}q((\emph{\textrm{Im }}F)^kX)\langle X\rangle^{-2\alpha}\big) \in S(1,\langle X\rangle^{-2\alpha}dX^2). \num
\end{equation}
\end{lemma}

\bigskip

\noindent
\textit{Proof of Lemma~\ref{lem3}}. It is sufficient to check that 
\begin{equation}\label{giu17}\inc
\nabla\big(\textrm{Re }q((\textrm{Im }F)^kX)\langle X\rangle^{-2\alpha}\big) \in S_{\Omega}\big(\langle X \rangle^{-\alpha},\langle X \rangle^{-2\alpha}dX^2\big), \num
\end{equation}
where $\Omega$ is a small open neighborhood of  $\textrm{supp }f'\big(\textrm{Re }q((\textrm{Im }F)^kX)\langle X\rangle^{-2\alpha}\big).$ We deduce from (\ref{giu1}), (\ref{giu00.1}) and the fact that $\textrm{Re }q((\textrm{Im }F)^kX)$ is a quadratic form that 
$$\textrm{Re }q\big((\textrm{Im }F)^kX\big) \sim \langle X \rangle^{2\alpha}$$
and
$$|\nabla\big(\textrm{Re }q((\textrm{Im }F)^kX)\big)| \lesssim \textrm{Re }q((\textrm{Im }F)^kX)^{1/2} \lesssim \langle X\rangle^{\alpha},$$
on $\Omega$. By noticing that $0<\alpha \leq 1$, $\langle X \rangle^{r} \in S(\langle X \rangle^r,\langle X \rangle^{-2}dX^2)$, for any $r \in \rr$, and that the function 
$\textrm{Re }q((\textrm{Im }F)^kX)$ is just a quadratic form, we directly deduce (\ref{giu17}) from the previous estimates and the Leibniz's rule, since
$$\textrm{Re }q((\textrm{Im }F)^kX) \in S_{\Omega}\big(\langle X \rangle^{2\alpha},\langle X \rangle^{-2\alpha}dX^2\big). \ \Box$$

\bigskip

\begin{lemma}\label{lem2.15}
For all $s \in \rr$ and $0 \leq j \leq m-2$, we have
$$\emph{\textrm{Re }}q\big((\emph{\textrm{Im }}F)^{m-j-1}X\big)^{s} \in S_{\Omega}\Big(\emph{\textrm{Re }}q\big((\emph{\textrm{Im }}F)^{m-j-1}X\big)^{s},
\emph{\textrm{Re }}q\big((\emph{\textrm{Im }}F)^{m-j-1}X\big)^{-1}dX^2\Big),$$
if $\Omega$ is any open set where
$$\emph{\textrm{Re }}q\big((\emph{\textrm{Im }}F)^{m-j-1} X\big) \gtrsim \langle X \rangle^{\frac{2(2m-2j-1)}{2m+1}}.$$  
\end{lemma}

\bigskip

\noindent
\textit{Proof of Lemma~\ref{lem2.15}}. Recalling that the symbol $\textrm{Re }q\big((\textrm{Im }F)^{m-j-1}X\big)$ is a non-negative quadratic form and that we have from (\ref{giu00.1}) that 
\begin{equation}\label{mar1}\inc
\big|\nabla \textrm{Re }q\big((\textrm{Im }F)^{m-j-1}X\big)\big| \lesssim \textrm{Re }q\big((\textrm{Im }F)^{m-j-1}X\big)^{\frac{1}{2}}, \num
\end{equation}
which implies that for all $s \in \rr$,
\inc\begin{align*}\label{na3}
& \ \frac{\Big|\nabla\Big(\textrm{Re }q\big((\textrm{Im }F)^{m-j-1} X\big)^{s}\Big)\Big|}{\textrm{Re }q\big((\textrm{Im }F)^{m-j-1} X\big)^{s}}
\lesssim \frac{\big|\nabla\textrm{Re }q\big((\textrm{Im }F)^{m-j-1} X\big)\big|}{\textrm{Re }q\big((\textrm{Im }F)^{m-j-1} X\big)} \num \\
\lesssim & \ \textrm{Re }q\big((\textrm{Im }F)^{m-j-1} X\big)^{-\frac{1}{2}},
\end{align*}
on $\Omega$, we notice that the result of Lemma~\ref{lem2.15} is therefore a straightforward consequence of the Leibniz's rule.~$\Box$

\bigskip

\begin{lemma}\label{lem2.1}
Consider the function $\Psi_j$ defined in \emph{(\ref{giu39})} then for any $0 \leq j \leq m-2$,
$$\Psi_j \in S_{\Omega}\Big(1,\emph{\textrm{Re }}q\big((\emph{\textrm{Im }}F)^{m-j-1}X\big)^{-\frac{2m-2j-3}{2m-2j-1}}dX^2\Big),$$
if $\Omega$ is any open set where
$$\emph{\textrm{Re }}q\big((\emph{\textrm{Im }}F)^{m-j-1} X\big) \gtrsim \langle X \rangle^{\frac{2(2m-2j-1)}{2m+1}},$$  
which implies in particular that
$$\Psi_j \in S_{\Omega}\big(1,\langle X \rangle^{-\frac{2(2m-2j-3)}{2m+1}}dX^2\big).$$
\end{lemma}
\bigskip

\noindent
\textit{Proof of Lemma~\ref{lem2.1}}. We first notice from (\ref{giu13}) and (\ref{giu39}) that
$$\textrm{Re }q\big((\textrm{Im }F)^{m-j-2} X\big) \sim \textrm{Re }q\big((\textrm{Im }F)^{m-j-1}X\big)^{\frac{2m-2j-3}{2m-2j-1}},$$
on $\Omega \cap \textrm{supp }\Psi_j'$. 
Since from (\ref{giu00.1}), 
\inc\begin{align*}\label{na1}
 \big|\nabla\textrm{Re }q\big((\textrm{Im }F)^{m-j-2} X\big)\big| \lesssim & \ \textrm{Re }q\big((\textrm{Im }F)^{m-j-2} X\big)^{\frac{1}{2}} \num \\
\lesssim & \  \textrm{Re }q\big((\textrm{Im }F)^{m-j-1}X\big)^{\frac{2m-2j-3}{2(2m-2j-1)}},
\end{align*}
on $\Omega \cap \textrm{supp }\Psi_j'$, we deduce that the quadratic symbol 
$\textrm{Re }q\big((\textrm{Im }F)^{m-j-2} X\big)$ belongs to the class
\begin{equation}\label{na2}\inc
S_{\Omega \cap \textrm{supp}\Psi_j'}\Big(\textrm{Re }q\big((\textrm{Im }F)^{m-j-1}X\big)^{\frac{2m-2j-3}{2m-2j-1}},\frac{dX^2}{\textrm{Re }q\big((\textrm{Im }F)^{m-j-1}X\big)^{\frac{2m-2j-3}{2m-2j-1}}}\Big),\num
\end{equation}
It follows from Lemma~\ref{lem2.15} that 
\begin{multline*}
\frac{\textrm{Re }q\big((\textrm{Im }F)^{m-j-2} X\big)}{\textrm{Re }q\big((\textrm{Im }F)^{m-j-1} X\big)^{\frac{2m-2j-3}{2m-2j-1}}}
\in S_{\Omega \cap \textrm{supp}\Psi_j'}\Big(1,\frac{dX^2}{\textrm{Re }q\big((\textrm{Im }F)^{m-j-1}X\big)^{\frac{2m-2j-3}{2m-2j-1}}}\Big),
\end{multline*}
which implies that 
$$\Psi_j \in S_{\Omega}\Big(1,\textrm{Re }q\big((\textrm{Im }F)^{m-j-1}X\big)^{-\frac{2m-2j-3}{2m-2j-1}}dX^2\Big).$$
This ends the proof of Lemma~\ref{lem2.1}.~$\Box$

\bigskip

\begin{lemma}\label{lem2.2}
Consider the function $W_j$ defined in \emph{(\ref{giu40})} then for any $1 \leq j \leq m-1$,
$$W_j \in S_{\Omega}\Big(1,\emph{\textrm{Re }}q\big((\emph{\textrm{Im }}F)^{m-j-1} X\big)^{-1}dX^2\Big),$$
if $\Omega$ is any open set where
$$\emph{\textrm{Re }}q\big((\emph{\textrm{Im }}F)^{m-j-1} X\big) \gtrsim \langle X \rangle^{\frac{2(2m-2j-1)}{2m+1}},$$  
which implies in particular that 
$$W_j \in S_{\Omega}\big(1,\langle X \rangle^{-\frac{2(2m-2j-1)}{2m+1}}dX^2\big).$$
\end{lemma}

\bigskip

\noindent
\textit{Proof of Lemma~\ref{lem2.2}}. By noticing from (\ref{giu15}) and (\ref{giu40}) that
$$\textrm{Re }q\big((\textrm{Im }F)^{m-j-1} X\big) \sim \textrm{Re }q\big((\textrm{Im }F)^{m-j}X\big)^{\frac{2m-2j-1}{2m-2j+1}}$$
and
$$\textrm{Re }q\big((\textrm{Im }F)^{m-j} X\big) \gtrsim \langle X \rangle^{\frac{2(2m-2j+1)}{2m+1}},$$
on $\Omega \cap \textrm{supp }W_j'$, and that the two derivatives $\psi'$ and $w_2'$ of the functions appearing in (\ref{giu39}) and (\ref{giu40}) have similar types of support as the function defined in (\ref{giu14}), we notice that we are exactly in the setting studied in Lemma~\ref{lem2.1} with $j$ replaced by $j-1$. We therefore deduce the result of Lemma~\ref{lem2.2} from our analysis led in the proof of Lemma~\ref{lem2.1}.~$\Box$

\bigskip

\begin{lemma}\label{lem2.21}
If $l_1$, $l_2 \in \nn$ then 
$$\big|\emph{\textrm{Re }}q\big((\emph{\textrm{Im }}F)^{l_1}X;(\emph{\textrm{Im }}F)^{l_2}X\big)\big| \leq \emph{\textrm{Re }}q\big((\emph{\textrm{Im }}F)^{l_1}X\big)^{\frac{1}{2}}
\emph{\textrm{Re }}q\big((\emph{\textrm{Im }}F)^{l_2}X\big)^{\frac{1}{2}},
$$
$$\big|\nabla \emph{\textrm{Re }}q\big((\emph{\textrm{Im }}F)^{l_1}X;(\emph{\textrm{Im }}F)^{l_2}X\big)\big| \lesssim \emph{\textrm{Re }}q\big((\emph{\textrm{Im }}F)^{l_1}X\big)^{\frac{1}{2}}+\emph{\textrm{Re }}q\big((\emph{\textrm{Im }}F)^{l_2}X\big)^{\frac{1}{2}}.
$$\end{lemma}

\bigskip

\noindent
\textit{Proof of Lemma~\ref{lem2.21}}. By reason of symmetry, we can assume in the following that $l_1 \leq l_2$. Recalling that the quadratic form $\textrm{Re }q$ is non-negative, the first estimate is a direct consequence of the Cauchy-Schwarz inequality.
About the second estimate, we recall from (\ref{lol10}) that the Hamilton map of the quadratic form
$$\textrm{Re }q\big((\textrm{Im }F)^{l_1}X;(\textrm{Im }F)^{l_2}X\big),$$ 
is 
$$\frac{1}{2}\big((-1)^{l_1}(\textrm{Im }F)^{l_1}\textrm{Re }F(\textrm{Im }F)^{l_2}+(-1)^{l_2}(\textrm{Im }F)^{l_2}\textrm{Re }F(\textrm{Im }F)^{l_1}\big).$$ 
A direct computation as in (3.18) of~\cite{mz} shows that
\inc\begin{multline*}\label{mari1}
\nabla \textrm{Re }q\big((\textrm{Im }F)^{l_1}X;(\textrm{Im }F)^{l_2}X\big)=(-1)^{l_1+1}\sigma(\textrm{Im }F)^{l_1}\textrm{Re }F(\textrm{Im }F)^{l_2}X \num\\
+(-1)^{l_2+1}\sigma(\textrm{Im }F)^{l_2}\textrm{Re }F(\textrm{Im }F)^{l_1}X
\end{multline*}
where
$$\sigma=\begin{pmatrix}
0 & I_n\\
-I_n & 0
\end{pmatrix}.
$$ 
The notation $I_n$ stands here for the $n$ by $n$ identity matrix. We deduce from (\ref{giu00.1}) and (\ref{mari1}) that for any $k \in \nn$,
\inc\begin{equation}\label{kee4}
|(\textrm{Im }F)^{k}\textrm{Re }F(\textrm{Im }F)^{k}X| \lesssim  \big|\nabla \textrm{Re }q\big((\textrm{Im }F )^{k}X\big)\big|
\lesssim   \textrm{Re }q\big((\textrm{Im }F )^{k}X\big)^{\frac{1}{2}}.\num
\end{equation}
By using twice the estimate (\ref{kee4}) with respectively $X$ and $(\textrm{Im }F)^{l_2-l_1}X$, and the index $k=l_1$, we deduce from (\ref{mari1}) the second estimate of Lemma~\ref{lem2.21}.~$\Box$

\bigskip

\begin{lemma}\label{lem2.3}
Consider the quadratic form $r_{m-j-1}$ defined in \emph{(\ref{giu11})} then for any $0 \leq j \leq m-2$,
$$\frac{r_{m-j-1}(X)}{\emph{\textrm{Re }}q\big((\emph{\textrm{Im }}F )^{m-j-1}X\big)^{\frac{2m-2j-2}{2m-2j-1}}} \in 
S_{\Omega}\Big(1,\emph{\textrm{Re }}q\big((\emph{\textrm{Im }}F )^{m-j-1}X\big)^{-\frac{2m-2j-3}{2m-2j-1}}dX^2\Big),$$
if $\Omega$ is any open set where
$$\emph{\textrm{Re }}q\big((\emph{\textrm{Im }}F)^{m-j-1} X\big) \gtrsim \langle X \rangle^{\frac{2(2m-2j-1)}{2m+1}}$$ 
and  
$$\emph{\textrm{Re }}q\big((\emph{\textrm{Im }}F)^{m-j-2} X\big) \lesssim \emph{\textrm{Re }}q\big((\emph{\textrm{Im }}F)^{m-j-1} X\big)^{\frac{2m-2j-3}{2m-2j-1}},$$  
which implies in particular that 
$$\frac{r_{m-j-1}(X)}{\emph{\textrm{Re }}q\big((\emph{\textrm{Im }}F )^{m-j-1}X\big)^{\frac{2m-2j-2}{2m-2j-1}}} \in S_{\Omega}\big(1,\langle X \rangle^{-\frac{2(2m-2j-3)}{2m+1}}dX^2\big).$$
\end{lemma}

\bigskip

\noindent
\textit{Proof of Lemma~\ref{lem2.3}}. Since from Lemma~\ref{lem2.21},
$$|r_{m-j-1}(X)| \lesssim \textrm{Re }q\big((\textrm{Im }F )^{m-j-1}X\big)^{\frac{2m-2j-2}{2m-2j-1}}$$
and
\begin{align*}
|\nabla r_{m-j-1}(X)| \lesssim & \ \textrm{Re }q\big((\textrm{Im }F )^{m-j-1}X\big)^{\frac{1}{2}}+\textrm{Re }q\big((\textrm{Im }F )^{m-j-2}X\big)^{\frac{1}{2}} \\
\lesssim & \  \textrm{Re }q\big((\textrm{Im }F )^{m-j-1}X\big)^{\frac{1}{2}},
\end{align*} 
on $\Omega$, we get that the quadratic form $r_{m-j-1} $ belongs to the symbol class
$$S_{\Omega}\Big(\textrm{Re }q\big((\textrm{Im }F )^{m-j-1}X\big)^{\frac{2m-2j-2}{2m-2j-1}},\textrm{Re }q\big((\textrm{Im }F )^{m-j-1}X\big)^{-\frac{2m-2j-3}{2m-2j-1}}dX^2\Big).$$
One can then deduce the result of Lemma~\ref{lem2.3} from Lemma~\ref{lem2.15}.~$\Box$

\bigskip

When adding a large parameter $\Lambda_j \geq 1$ in the description of the open set $\Omega$, a straightforward adaptation of the proof of the previous lemma gives the 
following $L^{\infty}(\Omega)$ estimate with respect to this parameter.

\bigskip

\begin{lemma}\label{lem2.331}
Consider the quadratic form $r_{m-j-1}$ defined in \emph{(\ref{giu11})} then for any $0 \leq j \leq m-2$,
$$\big\|\emph{\textrm{Re }}q\big((\emph{\textrm{Im }}F )^{m-j-1}X\big)^{-\frac{2m-2j-2}{2m-2j-1}}r_{m-j-1}(X)   \big\|_{L^{\infty}(\Omega)} \lesssim \Lambda_j^{-\frac{1}{2}},$$
if $\Omega$ is any open set where
$$\emph{\textrm{Re }}q\big((\emph{\textrm{Im }}F)^{m-j-1} X\big) \gtrsim \langle X \rangle^{\frac{2(2m-2j-1)}{2m+1}}$$ 
and  
$$\emph{\textrm{Re }}q\big((\emph{\textrm{Im }}F)^{m-j-2} X\big) \lesssim \Lambda_j^{-1} \emph{\textrm{Re }}q\big((\emph{\textrm{Im }}F)^{m-j-1} X\big)^{\frac{2m-2j-3}{2m-2j-1}}.$$  
\end{lemma}

\bigskip

In the following lemmas, we shall carefully study the dependence of the estimates with respect to the large parameter $\Lambda_j \geq 1$.

\bigskip

\begin{lemma}\label{lem2.4}
For any $0 \leq j \leq m-2$, we have for all $X \in \Omega$,
\begin{multline*}
\left|\frac{H_{\emph{\textrm{Im}}q}\ r_{m-j-1}(X)}{\emph{\textrm{Re }}q\big((\emph{\textrm{Im }}F )^{m-j-1}X\big)^{\frac{2m-2j-2}{2m-2j-1}}} 
-2\emph{\textrm{Re }}q\big((\emph{\textrm{Im }}F)^{m-j-1} X\big)^{\frac{1}{2m-2j-1}}\right| \\
\lesssim \Lambda_j^{-\frac{1}{2}} \emph{\textrm{Re }}q\big((\emph{\textrm{Im }}F)^{m-j-1} X\big)^{\frac{1}{2m-2j-1}},
\end{multline*}
if $\Omega$ is any open set where
$$\emph{\textrm{Re }}q\big((\emph{\textrm{Im }}F)^{m-j-1} X\big) \gtrsim \langle X \rangle^{\frac{2(2m-2j-1)}{2m+1}},$$
$$\emph{\textrm{Re }}q\big((\emph{\textrm{Im }}F)^{m-j-2} X\big) \lesssim \Lambda_j^{-1} \emph{\textrm{Re }}q\big((\emph{\textrm{Im }}F)^{m-j-1} X\big)^{\frac{2m-2j-3}{2m-2j-1}},$$
$$\emph{\textrm{Re }}q\big((\emph{\textrm{Im }}F)^{m-j} X\big) \lesssim \emph{\textrm{Re }}q\big((\emph{\textrm{Im }}F)^{m-j-1} X\big)^{\frac{2m-2j+1}{2m-2j-1}}.$$
\end{lemma}
\bigskip

\noindent
\textit{Proof of Lemma~\ref{lem2.4}}. We begin by writing from (\ref{giu11}) and Lemma~\ref{lem2} that 
\inc\begin{multline*} \label{kee1}
H_{\textrm{Im}q}\ r_{m-j-1}(X)=2\textrm{Re }q\big((\textrm{Im }F )^{m-j-1}X\big)\\
+2\textrm{Re }q\big((\textrm{Im }F )^{m-j}X;(\textrm{Im }F )^{m-j-2}X \big).\num
\end{multline*}
Lemma~\ref{lem2.4} is then a consequence of the following estimate 
\begin{align*}
& \ \big|\textrm{Re }q\big((\textrm{Im }F )^{m-j}X;(\textrm{Im }F )^{m-j-2}X \big)\big| \\
\leq & \ \textrm{Re }q\big((\textrm{Im }F )^{m-j}X\big)^{\frac{1}{2}}
\textrm{Re }q\big((\textrm{Im }F )^{m-j-2}X \big)^{\frac{1}{2}}\\
\lesssim & \ \Lambda_j^{-\frac{1}{2}}\textrm{Re }q\big((\textrm{Im }F )^{m-j-1}X\big),
\end{align*}
fulfilled on $\Omega$ that we obtain from Lemma~\ref{lem2.21}.~$\Box$

\bigskip

\begin{lemma}\label{lem2.41}
For any $0 \leq j \leq m-2$, we have for all $X \in \Omega$,
\begin{multline*}
\Big|\emph{\textrm{Re }}q\big((\emph{\textrm{Im }}F )^{m-j-1}X\big)^{\frac{2m-2j-2}{2m-2j-1}}H_{\emph{\textrm{Im}}q}\Big(\emph{\textrm{Re }}q\big((\emph{\textrm{Im }}F)^{m-j-1}X\big)^{-\frac{2m-2j-2}{2m-2j-1}} \Big)\Big| \\ 
\lesssim \emph{\textrm{Re }}q\big((\emph{\textrm{Im }}F)^{m-j-1} X\big)^{\frac{1}{2m-2j-1}},          
\end{multline*}
if $\Omega$ is any open set where
$$\emph{\textrm{Re }}q\big((\emph{\textrm{Im }}F)^{m-j-1} X\big) \gtrsim \langle X \rangle^{\frac{2(2m-2j-1)}{2m+1}},$$
$$\emph{\textrm{Re }}q\big((\emph{\textrm{Im }}F)^{m-j-2} X\big) \lesssim \Lambda_j^{-1}\emph{\textrm{Re }}q\big((\emph{\textrm{Im }}F)^{m-j-1} X\big)^{\frac{2m-2j-3}{2m-2j-1}},$$
$$\emph{\textrm{Re }}q\big((\emph{\textrm{Im }}F)^{m-j} X\big) \lesssim \emph{\textrm{Re }}q\big((\emph{\textrm{Im }}F)^{m-j-1} X\big)^{\frac{2m-2j+1}{2m-2j-1}}.$$
\end{lemma}
\bigskip

\noindent
\textit{Proof of Lemma~\ref{lem2.41}}. We begin by writing from Lemma~\ref{lem2} that 
\inc\begin{equation} \label{kee1.001}
H_{\textrm{Im}q}\ \textrm{Re }q\big((\textrm{Im }F)^{m-j-1}X\big)=4\textrm{Re }q\big((\textrm{Im }F )^{m-j-1}X;(\textrm{Im }F)^{m-j}X\big). \num
\end{equation}
Since 
\begin{align*}
& \ \textrm{Re }q\big((\textrm{Im }F )^{m-j-1}X\big)^{\frac{2m-2j-2}{2m-2j-1}}H_{\textrm{Im}q}\Big(\textrm{Re }q\big((\textrm{Im }F)^{m-j-1}X\big)^{-\frac{2m-2j-2}{2m-2j-1}} \Big)\\
= & \ -\frac{2m-2j-2}{2m-2j-1}\frac{H_{\textrm{Im}q}\ \textrm{Re }q\big((\textrm{Im }F )^{m-j-1}X\big)}{\textrm{Re }q\big((\textrm{Im }F )^{m-j-1}X\big)},
\end{align*}
Lemma~\ref{lem2.41} is then a consequence of the following estimate
\inc\begin{align*}\label{spea2}
& \ \num \big|\textrm{Re }q\big((\textrm{Im }F )^{m-j-1}X;(\textrm{Im }F )^{m-j}X \big)\big| \\
\leq & \ \textrm{Re }q\big((\textrm{Im }F )^{m-j-1}X\big)^{\frac{1}{2}}
\textrm{Re }q\big((\textrm{Im }F )^{m-j}X \big)^{\frac{1}{2}}\\
\lesssim & \ \textrm{Re }q\big((\textrm{Im }F )^{m-j-1}X\big)^{1+\frac{1}{2m-2j-1}}, 
\end{align*}
fulfilled on $\Omega$ that we obtain from Lemma~\ref{lem2.21}.~$\Box$

\bigskip

\begin{lemma}\label{lem2.42}
Consider the functions $\Psi_j$ and $W_{j+1}$ defined in \emph{(\ref{giu39})} and \emph{(\ref{giu40})} then for any $0 \leq j \leq m-2$, we have for all $X \in \Omega$,
$$|H_{\emph{\textrm{Im}}q}\Psi_j(X)| \lesssim \Lambda_j^{\frac{1}{2}}
\emph{\textrm{Re }}q\big((\emph{\textrm{Im }}F)^{m-j-1} X\big)^{\frac{1}{2m-2j-1}}W_{j+1}(X),$$
if $\Omega$ is any open set where
$$\emph{\textrm{Re }}q\big((\emph{\textrm{Im }}F)^{m-j-1} X\big) \gtrsim \langle X \rangle^{\frac{2(2m-2j-1)}{2m+1}},$$
$$\emph{\textrm{Re }}q\big((\emph{\textrm{Im }}F)^{m-j-2} X\big) \lesssim \Lambda_j^{-1}\emph{\textrm{Re }}q\big((\emph{\textrm{Im }}F)^{m-j-1} X\big)^{\frac{2m-2j-3}{2m-2j-1}},$$
$$\emph{\textrm{Re }}q\big((\emph{\textrm{Im }}F)^{m-j} X\big) \lesssim \emph{\textrm{Re }}q\big((\emph{\textrm{Im }}F)^{m-j-1} X\big)^{\frac{2m-2j+1}{2m-2j-1}}.$$
\end{lemma}
\bigskip

\noindent
\textit{Proof of Lemma~\ref{lem2.42}}. We begin by noticing from (\ref{giu40}) and (\ref{bell1}) that
\begin{equation}\label{spea10}\inc
\left|\psi'\left(\frac{\Lambda_j \textrm{Re }q\big((\textrm{Im }F)^{m-j-2}X\big)}{\textrm{Re }q\big((\textrm{Im }F)^{m-j-1}X\big)^{\frac{2m-2j-3}{2m-2j-1}}}\right)\right| \lesssim
W_{j+1}(X),\num
\end{equation}
and by writing from Lemma~\ref{lem2} that 
\inc\begin{equation} \label{kee1.002}
H_{\textrm{Im}q}\ \textrm{Re }q\big((\textrm{Im }F)^{m-j-2}X\big)=4\textrm{Re }q\big((\textrm{Im }F )^{m-j-2}X;(\textrm{Im }F)^{m-j-1}X\big). \num
\end{equation}
It follows from Lemma~\ref{lem2.21} that for all $X \in \Omega$,
\inc\begin{align*}\label{spea1}
& \ \num \big|\textrm{Re }q\big((\textrm{Im }F )^{m-j-2}X;(\textrm{Im }F )^{m-j-1}X \big)\big| \\
\leq & \ \textrm{Re }q\big((\textrm{Im }F )^{m-j-2}X\big)^{\frac{1}{2}}
\textrm{Re }q\big((\textrm{Im }F )^{m-j-1}X \big)^{\frac{1}{2}}\\
\lesssim & \ \Lambda_j^{-\frac{1}{2}}\textrm{Re }q\big((\textrm{Im }F )^{m-j-1}X\big)^{\frac{2m-2j-2}{2m-2j-1}}. 
\end{align*}
Then, by writing that 
\begin{multline*}
H_{\textrm{Im}q}\left(\frac{\Lambda_j \textrm{Re }q\big((\textrm{Im }F)^{m-j-2}X\big)}{\textrm{Re }q\big((\textrm{Im }F)^{m-j-1}X\big)^{\frac{2m-2j-3}{2m-2j-1}}}\right)=
\frac{\Lambda_j H_{\textrm{Im}q} \ \textrm{Re }q\big((\textrm{Im }F)^{m-j-2}X\big)}{\textrm{Re }q\big((\textrm{Im }F)^{m-j-1}X\big)^{\frac{2m-2j-3}{2m-2j-1}}}
\\ -\frac{2m-2j-3}{2m-2j-1}\frac{\Lambda_j\textrm{Re }q\big((\textrm{Im }F )^{m-j-2}X\big)H_{\textrm{Im}q}\ \textrm{Re }q\big((\textrm{Im }F )^{m-j-1}X\big)}{\textrm{Re }q\big((\textrm{Im }F )^{m-j-1}X\big)^{1+\frac{2m-2j-3}{2m-2j-1}}}.
\end{multline*}
Lemma~\ref{lem2.42} is a consequence of (\ref{giu39}), (\ref{kee1.001}), (\ref{spea2}), (\ref{kee1.002}), (\ref{spea1}) and (\ref{spea10}), since 
$$\textrm{Re }q\big((\textrm{Im }F)^{m-j-2} X\big) \sim \Lambda_j^{-1}\textrm{Re }q\big((\textrm{Im }F)^{m-j-1} X\big)^{\frac{2m-2j-3}{2m-2j-1}},$$
on the support of $\Psi_j'$.~$\Box$

\bigskip

\begin{lemma}\label{lem2.435}
For $m \geq 2$, consider the function $\tilde{W}_0$ defined in \emph{(\ref{giu41})} then for all $X \in \rr^{2n}$,
$$|H_{\emph{\textrm{Im}}q}  \tilde{W}_0(X)| \lesssim \langle X \rangle^{\frac{2}{2m+1}}.$$
\end{lemma}
\bigskip

\noindent
\textit{Proof of Lemma~\ref{lem2.435}}. 
Since $|\nabla \textrm{Im }q(X)| \lesssim \langle X \rangle$, because $\textrm{Im }q$ is a quadratic form, Lemma~\ref{lem2.435} is then a consequence of (\ref{giu15}), (\ref{giu41}) and Lemma~\ref{lem3}.~$\Box$

\bigskip

\begin{lemma}\label{lem2.43}
Consider the function $W_{j+1}$ defined in \emph{(\ref{giu40})} then for any $0 \leq j \leq m-2$, we have for all $X \in \Omega$, 
$$|H_{\emph{\textrm{Im}}q}W_{j+1}(X)| \lesssim \Lambda_j^{\frac{1}{2}}
\emph{\textrm{Re }}q\big((\emph{\textrm{Im }}F)^{m-j-1} X\big)^{\frac{1}{2m-2j-1}}\Psi_j(X),$$
if $\Omega$ is any open set where
$$\emph{\textrm{Re }}q\big((\emph{\textrm{Im }}F)^{m-j-1} X\big) \gtrsim \langle X \rangle^{\frac{2(2m-2j-1)}{2m+1}},$$
$$\emph{\textrm{Re }}q\big((\emph{\textrm{Im }}F)^{m-j-2} X\big) \lesssim \Lambda_j^{-1}\emph{\textrm{Re }}q\big((\emph{\textrm{Im }}F)^{m-j-1} X\big)^{\frac{2m-2j-3}{2m-2j-1}},$$
$$\emph{\textrm{Re }}q\big((\emph{\textrm{Im }}F)^{m-j} X\big) \lesssim \emph{\textrm{Re }}q\big((\emph{\textrm{Im }}F)^{m-j-1} X\big)^{\frac{2m-2j+1}{2m-2j-1}}.$$
\end{lemma}
\bigskip

\noindent
\textit{Proof of Lemma~\ref{lem2.43}}. 
One can notice from (\ref{giu13}), (\ref{giu15}), (\ref{giu39}), (\ref{giu40}) and (\ref{bell1}) that
\begin{equation}\label{spea10}\inc
\forall \ 0 \leq j \leq n-2, \ \left|w_2'\left(\frac{\Lambda_j \textrm{Re }q\big((\textrm{Im }F)^{m-j-2}X\big)}{\textrm{Re }q\big((\textrm{Im }F)^{m-j-1}X\big)^{\frac{2m-2j-3}{2m-2j-1}}}\right)\right| \lesssim
\Psi_{j}(X),\num
\end{equation}
and that the derivatives of  $\Psi_j$ and $W_{j+1}$ are exactly the same types of functions. It follows that Lemma~\ref{lem2.43} is just a straightforward consequence of Lemma~\ref{lem2.42}.~$\Box$

\bigskip
\bigskip

\noindent
\textsc{Department of Mathematics,
Imperial College London,
Huxley Building, 180 Queen's Gate, 
London SW7 2AZ, UK}\\
\textit{E-mail address:} \textbf{k.pravda-starov@imperial.ac.uk}

\end{document}